\documentclass[a4paper]{amsart}
\usepackage{microtype}%if unwanted, comment out or use option "draft"
\usepackage{calc}
\usepackage{enumerate}
\usepackage{ifthen}
\usepackage{epsfig}
\usepackage{graphicx}
\usepackage{latexsym}
\usepackage{amssymb,amsmath}

\begin{document}
%%%%%%%%%%%%%%%%%%%%%%%%%%%%%%%%%%%%%%%%%%%%%%%%
%\title{Congruences for restricted partition functions} %TODO Please add
%\title{Counting labeled structures  \newline with hard-wired constants} %TODO Please add
\title{MC-finiteness of restricted set partition functions}
%------------------------------------------
\author{Y. Filmus}
\address{Faculty of Computer Science
Technion-Israel Institute of Technology, Haifa Israel}
\email{yuvalfi@technion.ac.il}
%------------------------------------------
\author{E. Fischer}
\address{Faculty of Computer Science
Technion-Israel Institute of Technology, Haifa Israel}
\email{eldar@cs.technion.ac.il}
%------------------------------------------
\author{J.A. Makowsky}
\address{Faculty of Computer Science
Technion-Israel Institute of Technology, Haifa Israel}
\email{janos@cs.technion.ac.il}
%------------------------------------------
\author{V. Rakita}
\address{Faculty of Mathematics
Technion-Israel Institute of Technology, Haifa Israel}
\email{vsevolod@campus.technion.ac.il}
%--------------------------------------------------

%\subjclass{tbc}
% mandatory: Please choose ACM 1998 classifications from http://www.acm.org/about/class/ccs98-html . 
%E.g., cite as "F.1.1 Models of Computation". 
% mandatory: Please provide 1-5 keywords
%%%%%%%%%%%%%%%%%%%%%%%%%%%%%%%%%%%%%%%%%%%%%%%%%%%%%%%%%

%Editor-only macros (do not touch as author)
%%%%%%%%%%%%%%%%%%%%%%%%%%%%%%%%%%%
%\serieslogo{}%please provide filename (without suffix)
%\volumeinfo%(easychair interface)
  %{Billy Editor, Bill Editors}% editors
  %{2}% number of editors: 1, 2, ....
  %{TTL-2015, Tools for Teaching Logic }% event
  %{1}% volume
  %{1}% issue
  %{1}% starting page number
%\EventShortName{}
%\DOI{10.4230/LIPIcs.xxx.yyy.p}% to be completed by the volume editor
%%%%%%%%%%%%%%%%%%%%%%%%%%%%%%%%%%%%%%%%%%%%%%%%%%%%%%%%%
%--------------- end lipic----------------------------
%\usepackage{theorem}
%-------------------------

\newcommand{\angl}[1]{\left\langle #1 \right\rangle}
\newcommand{\card}[3]{card_{\mathcal{#1},\bar{#2}}(#3(\bar{#2}))}
\newif\ifmargin
%\margintrue
\marginfalse
\newif\ifshort
\shorttrue
\newif\ifskip
\skiptrue
\newcommand{\NN}{\mathbb{N}}
\newcommand{\ZZ}{\mathbb{Z}}
%-----------------------------------------------------------------------------
%\input{rj-macros}
\newtheorem{theorem}{Theorem}
\newtheorem{proposition}[theorem]{\bf Proposition}
\newtheorem{examples}[theorem]{\bf Examples}
\newtheorem{example}[theorem]{\bf Example}
\newtheorem{problem}{\bf Problem}
\newtheorem{remark}[theorem]{\bf Remark}
\newtheorem{remarks}[theorem]{\bf Remarks}
\newtheorem{definition}{Definition}
\newtheorem{corollary}[theorem]{Corollary}
%---------------------------------------------------
\newtheorem{lesson}{Lesson}
\newtheorem{defi}{Definition}[section]
\newtheorem{conjecture}{Conjecture}
\newtheorem{ex}{Example}[section]
\newtheorem{lemma}[theorem]{Lemma}
\newtheorem{coro}[theorem]{Corollary}
\newtheorem{conj}[theorem]{Conjecture}
\newtheorem{cons}[theorem]{Consequence}
\newtheorem{obs}[theorem]{Observation}
\newtheorem{claim}[theorem]{Claim}
\newtheorem{fact}[theorem]{Fact}
\newtheorem{oproblem}[theorem]{Problem}
% CUSTOMIZING NUMBERED LISTS
\newenvironment{renumerate}{\begin{enumerate}}{\end{enumerate}}
\renewcommand{\theenumi}{\roman{enumi}}
\renewcommand{\labelenumi}{(\roman{enumi})}
\renewcommand{\labelenumii}{(\roman{enumi}.\alph{enumii})}
%----------------------------------------------------------
\renewcommand{\tilde}{\widetilde}
\renewcommand{\bar}{\overline}
%MATH MODE ABBREVIATIONS
\newcommand{\dd}{\mathrm{D}}
\newcommand{\WFF}{\mathrm{WFF}}
\newcommand{\SOL}{\mathrm{SOL}}
\newcommand{\FOL}{\mathrm{FOL}}
\newcommand{\MSOL}{\mathrm{MSOL}}
\newcommand{\CMSOL}{\mathrm{CMSOL}}
\newcommand{\CFOL}{\mathrm{CFOL}}
\newcommand{\IFPL}{\mathrm{IFPL}}
\newcommand{\FPL}{\mathrm{FPL}}
\newcommand{\SEN}{\mbox{\bf SEN}}
\newcommand{\WFTF}{\mbox{\bf WFTF}}
\newcommand{\TFOF}{\mbox{\bf TFOF}}
\newcommand{\TFOL}{\mbox{\bf TFOL}}
\newcommand{\FOF}{\mbox{\bf FOF}}
\newcommand{\NNF}{\mbox{\bf NNF}}
\newcommand{\N}{{\mathbb N}}
\newcommand{\bN}{{\mathbb N}}
\newcommand{\bR}{{\mathbb R}}
\newcommand{\HF}{\mbox{\bf HF}}
\newcommand{\CNF}{\mbox{\bf CNF}}
\newcommand{\PNF}{\mbox{\bf PNF}}
\newcommand{\QF}{\mbox{\bf QF}}
\newcommand{\DNF}{\mbox{\bf DNF}}
\newcommand{\DISJ}{\mbox{\bf DISJ}}
\newcommand{\CONJ}{\mbox{\bf CONJ}}
\newcommand{\Ass}{\mbox{Ass}}
\newcommand{\Var}{\mbox{Var}}
\newcommand{\Support}{\mbox{Support}}
\newcommand{\V}{\mbox{\bf Var}}
\newcommand{\fA}{{\mathfrak A}}
\newcommand{\fB}{{\mathfrak B}}
\newcommand{\fN}{{\mathfrak N}}
\newcommand{\fZ}{{\mathfrak Z}}
\newcommand{\fQ}{{\mathfrak Q}}
\newcommand{\Aa}{{\mathfrak A}}
\newcommand{\Bb}{{\mathfrak B}}
\newcommand{\Cc}{{\mathfrak C}}
\newcommand{\Gg}{{\mathfrak G}}
\newcommand{\Ww}{{\mathfrak W}}
\newcommand{\Rr}{{\mathfrak R}}
\newcommand{\Nn}{{\mathfrak N}}
\newcommand{\Zz}{{\mathfrak Z}}
\newcommand{\Qq}{{\mathfrak Q}}
\newcommand{\F}{{\mathbf F}}
\newcommand{\T}{{\mathbf T}}
\newcommand{\Z}{{\mathbb Z}}
\newcommand{\R}{{\mathbb R}}
\newcommand{\C}{{\mathbb C}}
\newcommand{\Q}{{\mathbb Q}}
\newcommand{\bP}{{\mathbf P}}
\newcommand{\bPH}{{\mathbf{PH}}}
\newcommand{\bNP}{{\mathbf{NP}}}
\newcommand{\MT}{\mbox{MT}}
\newcommand{\TT}{\mbox{TT}}
\newcommand{\cL}{\mathcal{L}}
%-----------------------------------------------------------------------------
\begin{abstract}
A sequence $s(n)$ of integers is MC-finite if for every $m \in \N$ the sequence $s^m(n) = s(n) \bmod{m}$
is ultimately periodic. We discuss various ways of proving and disproving MC-finiteness.
Our examples are mostly taken from set partition functions, but our methods can be applied to 
many more integer sequences.
%--------------------------------------------------
\end{abstract}
\maketitle
%\today
%-----------------------------------------------------------------
\sloppy
%---------------------------------------------
%\input{rj-files}
\footnotesize
\today
\tableofcontents
\normalsize
\newpage
%-------------------------------------------------------------------
%\input{rj-intro}
\section{Introduction}
\label{se:intro}
\ifmargin
\marginpar{s-intro}
\else\fi

\subsection{Goal of this paper}
Given a sequence of integers $s(n)$ with some combinatorial interpretation, one wonders what can be said
about the  sequence $s(n)$. Ideally, we would like to have an explicit formula for $s(n)$, 
or some recurrence relation with coefficients being constant or polynomial in $n$.
Second best is an asymptotic description of $s(n)$. 

%The monograph \cite{bk:HararyPalmer} is a rich source  for sequences which  enumerate objects in graphs of order $n$.
%If neither of these descriptions are known about $s(n)$, 
We could instead look at the
sequence $s^m(n) \equiv s(n) \bmod{m}$ and try to describe $s^m(n)$.
If for every  modulus $m$ the sequence $s^m(n)$ is ultimately periodic, we say that $s(m)$
is {\em MC-finite}. We consider MC-finiteness a legitimate topic in the study of integer sequences.
MC-finiteness appears under this name only since the publication of \cite{makowsky2010application} in 2010.
Without its name, the concept appears in the literature before, but rarely, e.g., 
under the name of {\em supercongruence}
\cite{banderier2017right,banderier2019period}.
The four substantial monographs
on integer sequences published after 2000 
do not mention the concept at all, 
see \cite{everest2003recurrence,mansour2012combinatorics,mansour2012combinatorics,mezo2019combinatorics}.

All the sequences we discuss in this paper appear in 
{\em The On-Line Encyclopedia of Integer Sequences, OEIS, https://oeis.org/}, \cite{oeis}, with a number starting with $A$. 
We give these numbers with the first mention of the sequence, and list them also at the end of the paper.
Needless to say, our methods also apply to many other entries in OEIS.

This paper grew out of our attempts to show that the sequence $B_r(n)$ 
 of restricted Bell numbers
(only listed in OEIS for 
$r=2, A005493$ and
$r=3, A005494$)
and $S_r(n,k)$ of restricted Stirling numbers of the second kind $A143494-A143496$ 
introduced in \cite{broder1984r} are MC-finite.

The purpose of this paper is two-fold.
Its first part is mostly expository and written with the intent to popularize the study of MC-finiteness
for researchers in Integer Sequences. 
However, the statements that the examples chosen are MC-finite have not, to the best of our knowledge,
been stated before in the literature.
We have chosen our examples
in order to familiarize the reader with the two general methods to establish MC-finiteness.
The first is {\em logical methods}, pioneered by C. Blatter and E. Specker, 
\cite{blatter1981nombre,specker1990application,pr:BlatterSpecker84}, and further 
developed by two of the authors of this paper (EF and JAM), \cite{pr:FischerMakowsky03,ar:FischerMakowsky2022}.
%in understanding the modular behavior of integer sequences.
The second is a {\em combinatorial method} to prove MC-finiteness, also first suggested by E. Specker in
\cite{specker1990application}, and later independently by G. S\'enizergues \cite{senizergues2007sequences},
but only made precise in \cite{cadilhac2021polynomial}.
This method is based on the existence of finitely many mutual polynomial recurrence relations over $\Z$
used to define the integer sequence.
In a separate paper, these methods are applied to infinitely many integer sequences arising from
finite topologies \cite{topologies}. 

In this paper we investigate  MC-finiteness and counterexamples thereof of
%linear recurrence relations with constant coefficients and congruences of
integer sequences derived from counting various unrestricted and restricted set partitions and transitive relations.
Among the unrestricted cases we look at the Bell numbers $B(n)$, $A000110$, and the Stirling numbers
of the second kind $S(n, k_0)$, $A000453$. We also discuss
the number of linear quasi-orders (pre-orders) $LQ(n)$, $A000670$,
the number of quasi-orders (pre-orders) $Q(n)$, $A000798$,
the number of partial orders $P(n)$, $A001035$, and 
the number of transitive relations $T(n)$, $A006905$, on the set $[n]$.
The numbers $LQ(n)$ are called {\em ordered Bell numbers} or {\em Fubini numbers}, often denoted 
in the literature by $a(n)$ and also by $F(n)$.
For the unrestricted cases the results are seemingly new, or at least have not been stated before, 
but are simple consequences
of growth arguments and the logical method due to C. Blatter and E. Specker 
\cite{pr:BlatterSpecker84,specker1990application}, the {\em Specker-Blatter Theorem}.

Typical restricted cases, first introduced by A. Broder \cite{broder1984r} and further studied in \cite{benyi2019restricted}, 
are the Stirling numbers of the second kind $S_{A,r}(n,k)$,
which count the partitions of $[n+r]$ into $k+r$ blocks such that
the elements $i \leq r$ are all in different blocks and the size of each block is in $A \subseteq \NN$.
For $r=2$ see $A143494$. 
The Bell numbers $B_{A,r}(n)$ are defined as $\sum_k S_{A,r}(n,k)$, see $A005493$ for $r=2$ and $A005494$ for $r=3$.
The same restrictions can also be imposed on Stirling numbers of the second kind $S_{A,r}(n,k)$, 
and on all the unrestricted cases above.
For the restricted cases, the results are new and require non-trivial extensions of the Specker-Blatter
Theorem. 
The Catalan numbers $A000108$ also have an interpretation as set partitions. They count the number of
non-overlapping partitions, see \cite[Theorem 9.4]{roman2015introduction} and \cite[Chapter 10]{koshy2008catalan}.
Although this can be viewed as a restricted  version of the Bell numbers, our results do not apply to this case,
as we shall explain later.
%which are the main technical contribution of this paper.

%-------------------------------------------------------------------
%\input{rj-outline}
\subsection{Outline of the paper}
\ifmargin
\marginpar{s-outline}
\else\fi

In Section \ref{se:mcfinite} we introduce C-finiteness and its modular variant MC-finiteness.
In Section \ref{se:howtoprove} we discuss the methods for proving and disproving
C-finiteness and
MC-finiteness, and in Section \ref{se:immediate}
we present immediate consequences of the logical method for set partitions without 
positional restrictions and without restrictions on size of the blocks.
The three sections have tutorial character, although the MC-finiteness of the examples has not been
stated before in the literature.
In Sections \ref{se:restricted} and \ref{se:rproofs} we discuss  set partitions with
positional restrictions and restrictions on size of the blocks, and how new logical tools
are used to obtain 
C-finiteness and 
MC-finiteness in these cases.
We conclude the main part of the paper with Section \ref{se:conclu}, where
we present our conclusions and suggestions for further research, and in Section \ref{se:oeis} we list
the numbers of the discussed OEIS-sequences.
There are four appendices.
In Appendix \ref{ap:mc-finite} we discuss larger classes of polynomial recursive sequences and weaker versions of MC-finiteness.
In Appendix \ref{se:constants} we prove a special case of the main theorem from \cite{ar:FischerMakowsky2022}
which suffices for our results in Section \ref{se:rproofs}.
In Appendix \ref{se:c-finite} we give the details for proving C-finiteness of restricted Stirling numbers of the second kind.
Finally, in Appendix \ref{se:explicit}, we give an explicit computation of $S_A(n,k)$.

\ifskip\else
%-----------------------------------------------
In Section \ref{se:periodic}
we first deal with  the applications of the extension of the Specker-Blatter Theorem
to $\CMSOL$, Monadic Second Order Logic with modular counting. 
This allows us to show MC-finiteness 
for $B_{A,0}(n)$ 
and its variations 
for $A \subseteq \N$
ultimately periodic and
without hard-wired constants. 

In Section \ref{se:growth} we list some properties of restricted Stirling and Lah numbers
and give lower bounds for their growth. This shows that $B_{A,r}(n)$ and the corresponding Lah-numbers
are not C-finite provided that $A$ is infinite and ultimately periodic. 
Here a fixed finite number of hard-wired constants is allowed.

In Section \ref{se:fm} we give the background needed to prove C-finiteness using model theoretic methods
as described in \cite{fischer2008linear}.
This establishes C-finiteness for the restricted Stirling numbers of the second kind.

%In Section \ref{se:sb} we give the background needed to prove MC-finiteness using the methods of C. Blatter
%and E. Specker, \cite{pr:BlatterSpecker84,specker1990application} and its extensions from \cite{fischer2011application}.

In Section \ref{se:hard} we discuss hard-wired constants and
state a new extension of the Specker-Blatter Theorem, Theorem \cite{ar:FischerMakowsky2022} which shows how to handle
the case where a fixed finite set of hard-wired constants are allowed. The special case for 
Bell, Stirling and Lah numbers is proved in In Section \ref{se:constants}. 

% In Section \ref{se:mc} we prove all the cases where MC-finiteness holds.

Finally, in Section \ref{se:conclu}, we 
draw conclusions and suggest further lines of research concerning the complexity of computing
$S_{\phi}(n)$ and
$S^m_{\phi}(n)$, Problems 
\ref{problem-1},
\ref{problem-2} and
\ref{problem-3}.

There are two appendices.
In Appendix \ref{se:c-finite} we give a detailed prof of Theorem \ref{th:FM} and its applications.
In Appendix \ref{se:explicit} we give an explicit method of computation of $S_A(n,k_0)$ for fixed $k_0$ and
computable $A$.
\fi %skip

%-------------------------------------------------------------------
%\input{rj-methods}
\section{C-finite and MC-finite sequences of integers}
\label{se:mcfinite}
A sequence of integers $s(n)$ is {\em C-finite}\footnote{
These are also called constant-recursive sequences
or  linear-recursive sequences in the literature.
}
if there are constants $p, q \in \N$ and $c_i \in \Z, 0 \leq i \leq p-1$ such that for all $n \geq q$ the
linear recurrence relation
$$
s(n+p) = \sum_{i=0}^{p-1} c_i s(n+i), n \geq q,
$$
holds for $s(n)$.
C-finite sequences have limited growth, see e.g. \cite{everest2003recurrence,kauers2011concrete}:
\begin{proposition}
\label{prop:c-finite}
Let $s_n$ be a C-finite sequence of integers. Then there is $c \in \N^+$ such that for all $n \in \N$,
$a_n \leq 2^{cn}$.
\end{proposition}
Actually, a lot more can be said, see \cite{flajolet2009analytic}, but we do not need it for our purposes.

To prove that a sequence $s(n)$ of integers is not C-finite, we can use Proposition \ref{prop:c-finite}.
To prove that a sequence $s(n)$ of integers is C-finite, there are several methods:
One can try to find an explicit recurrence relation, one can exhibit a rational generating function,
or one can use a method based on model theory as described in 
\cite{fischer2008linear,fischer2011application}.
The last method will be briefly discussed in
Section \ref{se:fm} and further explained in Appendix \ref{se:c-finite}. It is referred to as method FM.

A sequence of integers $s(n)$ is {modular C-finite}, abbreviated as {\em MC-finite}, 
if for every $m \in \N$ there are constants $p_m, q_m \in \N^+$ such that
for every  $n \geq q_m$ there is a linear recurrence relation
$$
s(n+p_m) \equiv \sum_{i=0}^{p_m-1} c_{i,m}  s(n+i) \bmod{m}
$$
with constant coefficients $c_{i,m} \in \Z$. 
%and $n \geq n_0$.
Note that the coefficients $c_{i,m}$ and both $p_m$ and $q_m$  generally do depend on $m$.

We denote by $s^m(n)$ the sequence $s(n) \bmod{m}$. 
\begin{proposition}
The sequence $s(n)$o is MC-finite iff $s^m(n)$ is ultimately periodic for every $m$.
\end{proposition}
\begin{proof}
MC-finiteness implies periodicity. The converse is from \cite{reeds1985shift}.
\end{proof}

Clearly, if a sequence $s(n)$ is C-finite it is also MC-finite with $r_m=r$ and $c_{i,m}=c_i$ for all $m$.
The converse is not true, there are uncountably many MC-finite sequences, but only
countably many C-finite sequences with integer coefficients, see Proposition \ref{pr:many} below.
%Here are some typical examples:
\begin{examples}\ 
\label{ex:mc}
\begin{enumerate}[(i)]
\item
The Fibonacci sequence is C-finite.
\item
If $s(n)$ is C-finite it has at most simple exponential growth, by Proposition \ref{prop:c-finite}.
\item
The Bell numbers $B(n)$ are {\em not C-finite}, but are {\em MC-finite}.
\item
Let $f(n)$ be any integer sequence. The sequence $s_1(n)=2\cdot f(n)$ is ultimately periodic modulo $2$,
but not necessarily MC-finite.
\item
Let $g(n)$ be any integer sequence.
%grow arbitrarily fast. 
The sequence  $s_2(n) = n!\cdot g(n)$ is MC-finite.
\label{many-mc}
%We conclude that there are uncountably many monotonously increasing sequences which are MC-finite.
\item
The sequence $s_3(n)= \frac{1}{2} {2n \choose n}$ is not MC-finite: 
$s_3(n)$ is odd iff $n$ is a power of $2$, and otherwise it is even (Lucas, 1878).
A proof may be found in \cite[Exercise 5.61]{graham1989concrete} or in \cite{specker1990application}.
\item
The Catalan numbers $C(n) = \frac{1}{n+1}{2n \choose n}$ are not MC-finite,
since $C(n)$ is odd iff $n$ is a Mersenne number, i.e.,  $n = 2^m-1$ for some $m$,
see \cite[Chapter 13]{koshy2008catalan}.
%For a recent short proof of this, see \cite{KoshySalmassi}, and for an equivalent result from 1973, 
%see \cite[Theorem2]{alter1973binary}.
\item
\label{many-nonmc}
Let $p$ be a prime and $f(n)$ monotone increasing.
The sequence 
$$
s(n) = \begin{cases}
p^{f(n)} & n \neq p^{f(n)} \\
p^{f(n)}+1 & n = p^{f(n)}
\end{cases}
$$
is monotone increasing but  not ultimately periodic modulo $p$, hence not MC-finite.
%We conclude that there are uncountably many monotonously increasing sequences which are not MC-finite.
\end{enumerate}
\end{examples}

\begin{proposition}
\label{pr:many}
\begin{enumerate}[(i)]
\item
There are uncountably many monotone increasing sequences which are MC-finite, and uncountably many
which are not MC-finite.
\item
Almost all integer sequences are not MC-finite.
\end{enumerate}
\end{proposition}
\begin{proof}
(i) follows from
Examples \ref{ex:mc}
(\ref{many-mc}) and (\ref{many-nonmc}).
(ii) is shown in
Proposition \ref{pr:normal} in  Appendix \ref{se:ProofMC}.
\end{proof}

Although we are mostly interested in MC-finite sequences $s(n)$, it is natural to check in each example
whether the sequence $s(n)$ is also C-finite. In most examples the answer is negative.
However, Theorem \ref{th:c-finite} shows that for restricted Stirling numbers of the second kind
are all C-finite. We show this via a general method, Theorem \ref{th:FM}, without exhibiting a generating function
like in the classical case for $S(n,k)$.

\section{How to prove and disprove MC-finiteness}
\label{se:howtoprove}

\subsection{Polynomial recurrence relations}

In his 1988 paper \cite[Page 144]{specker1990application}, E. Specker notes the following:

\begin{quote}
In many known cases, [MC-finiteness] is a consequence of
polynomial recurrence relations 
$$f(n) = \sum_{i=1}^d P_i(n) f(n-i)$$
where $P_i$ are polynomials in $\Z[x]$.
\end{quote}
For $f(n) = n!$ this is obvious.

\begin{definition}
\begin{enumerate}[(i)]
\item
An integer sequence $s(n)$ is {\em holonomic over $\Z$} 
if there exist
polynomials $P_i \in \Z[x]$ with $P_1, P_k \neq 0$ such that
$$
s(n) = \sum_{i=1}^k P_i(n)s(n-i)
$$
\item
An integer sequence $s(n)$ is {\em polynomially recursive (PRS) over $\Z$} if
there exist $k \in \N$ integer sequences $s_i(n), 1 \leq i \leq k$ 
with $s(n) = s_1(n)$
and polynomials $P_i \in \Z[x_1, \ldots , x_k]$
such that  the following mutual recursion holds:
$$
s_i(n+1) = P_i(s_1(n), \ldots , s_k(n)), i = 1, \ldots k
$$
\item
An integer sequence $s(n)$ is {\em PRS  over $\Z$ and $n$} if
the polynomials also involve $n$ as an additional variable.
In other words $P_i \in \Z[x_1, \ldots , x_k, y]$ and
$$
s_i(n+1) = P_i(s_1(n), \ldots , s_k(n), n), i = 1, \ldots k
$$
\end{enumerate}
Actually, (ii) and (iii) are equivalent.
\end{definition}

We note that,
if $s(n)$ is an integer sequence which is polynomially recursive over $\Z$ and $n$ then $s(n)$
is holonomic over $\Z$.

In fact, the following is true:
\begin{theorem}
\label{th:MC}
If $s(n)$ is an integer sequence which is polynomially recursive over $\Z$ and $n$ then $s(n)$
is MC-finite. 
In particular, this is true also for integer sequences $s(n)$ 
holonomic over $\Z$.
\end{theorem}
The proof is given in Appendix \ref{se:ProofMC}. There we also briefly discuss weaker properties than
MC-finiteness, where the modular recurrence holds only for almost all $m \in \N^+$.

\begin{remarks}
\begin{enumerate}[(i)]
\item
In general, holonomic sequences are defined over  fields $\mathbb{F}$ rather than the ring $\Z$.
A good reference is \cite[Chapter 7]{kauers2011concrete}.
A theorem related to Theorem \ref{th:MC} for holonomic sequences can be found in \cite[Theorem 7]{banderier2017right},
see also \cite{banderier2019period}.
\item
In \cite{cadilhac2021polynomial}, polynomially recursive sequences are defined for rational numbers rather than
integers, and the polynomials are in $\Q[x_1, \ldots , x_k]$.
\end{enumerate}
\end{remarks}

The following examples, besides (v), are from \cite{cadilhac2021polynomial}.
\begin{examples}
\begin{enumerate}[(i)]
\item
The sequences $a(n) = n!$ with $a(n) = n\cdot a(n-1)$ and $a(0)=1$ is holonomic over $\Z$.
It is obviously MC-finite.
\item
The sequence
$a(n)= 2^{2^n}$ is polynomially recursive with $a(0)=2$ and $a(n) = a(n-1)^2$. 
It is not holonomic, since every holonomic sequence $a(n)$ is bounded by some $2^{p(n)}$
for some polynomial $p(n)$, see \cite{gerhold2004some}. 
It is easy to see that it is MC-finite, but it is also MC-finite by the Specker-Blatter Theorem below,
as it counts the number of ways one can interpret a unary predicate on $[n]$.
\item
The Catalan numbers $C_n$ are holonomic: $ (n+2)C_{n+1} = (4n+2)C_n $.
They are not holonomic over $\Z$, since they are not MC-finite. 
%see \ref{xxx}.
Furthermore, they are not polynomially recursive even if we allow rational numbers.
\item
The sequence $n^n$ is not polynomially recursive, but
it is MC-finite by the Specker-Blatter Theorem below.
\item
We show in Appendix \ref{ap:mc-finite} that
the sequence $A086714$ given by $a(0)=4, a(n+1) = {a(n) \choose 2}$ is not MC-finite but 
periodic modulo every odd number.
\end{enumerate}
\end{examples}

MC-finite sequences are closed under various arithmetic operations.
%\marginpar{check the last item}
\begin{proposition}
\label{mc-closure}
Let $a(n), b(n)$ be MC-finite sequences and $c \in \Z$.
\begin{enumerate}[(i)]
\item
Then $c\cdot a(n), a(n) + b(n), a(n) \cdot b(n)$ are MC-finite.
\item
If additionally, $b(n) \in \N^+$ and tends to infinity, 
$a(n)^{b(n)}$ is also MC-finite.
\item
Let  $A \subseteq \N^+$ be non-periodic and  $a(n) =2$ be a constant, hence MC-finite, sequence.
The sequence
$$
b(n) =
\begin{cases}
1 & n \in A \\
n!+1 & n \not \in A
\end{cases}
$$
is MC-finite and oscillates. However $a(n)^{b(n)}$ is not MC-finite. 
%Check for $n \bmod{4}$.
\end{enumerate}
\end{proposition}

%\begin{proof}[of Theorem \ref{th:MC}]
%See appendix
%\end{proof}

\subsection{A definability criterion}

In order to prove that a sequence $s(n)$ is MC-finite one can also use a method
due to E. Specker and C. Blatter from 1981
\cite{blatter1981nombre,pr:BlatterSpecker84,specker1990application}. 
It uses logical definability as a sufficient condition.
We denote by $\FOL$ first order logic, by $\MSOL$ monadic second order logic, and
by $\CMSOL$ the logic $\MSOL$ 
augmented with modular counting quantifiers. Details on the definition of $\CMSOL$ are given in Section
\ref{se:periodic}.
In its simplest form, the Specker Blatter Theorem can be stated as follows:

\begin{theorem}[Specker-Blatter Theorem]
\label{th:BS-CMSOL}
Let $S_{\phi}(n)$ be
the number of binary relations $R$ on a set $[n]$ which 
%can be defined in First Order Logic $\FOL$ by a formula $\phi$. 
satisfy a given formula $\phi \in \CMSOL$.
$S_{\phi}(n)$ is MC-finite, or equivalently, $S^m_{\phi}(n)$ is ultimately periodic for every $m$.
\end{theorem}
The original Specker-Blatter Theorem was stated for classes of structures with a finite set of binary relations
definable in Monadic Second Order Logic $\MSOL$. It also works with unary relations added. 
The extension to $\CMSOL$ is due to \cite{pr:FischerMakowsky03}.
%but not when constants are allowed. 
This method 
%will be described in Section \ref{se:bs} 
is abbreviated in the sequel by SB.

\subsection{Comparing the methods}
If one proves MC-finiteness for an integer sequence directly, the proof may be sometimes straightforward,
but also sometimes tricky, and not applicable to other sequences.
In contrast to this, Theorems \ref{th:MC} and \ref{th:BS-CMSOL} are meta-theorems. They only require to
check for some structural data about the sequence $s(n)$, recurrence relations or logical definability.
However, these meta-theorems  are only existence theorems, 
without explicitly giving the required coefficients $c_{i,m}$ which show MC-finiteness.

\begin{examples}
We note that the two meta-theorems cannot always be applied to the same integer sequences.
\begin{enumerate}[(i)]
\item
The sequence $s(n)=n^n$ counts the number of unary functions (as binary relations) from $[n]$ to $[n]$,
which is $\FOL$-definable, but it is not polynomially recursive, as shown in \cite{cadilhac2021polynomial}.
However, MC-finiteness can also be established directly without much effort.
\item
There are polynomially recursive sequences over $\Z$ (hence MC-finite) which grow as fast as $2^{2^n}$,
e.g., the sequence $a(0) = 2, a(n+1) = a(n)^2$ satisfies $a(n) = 2^{2^n}$. 
However, counting the number of $k$ binary relations on $[n]$ is bounded by $2^{kn^2}$.
Hence, Theorem \ref{th:BS-CMSOL} cannot be applied.
Again, MC-finiteness can also be established directly without much effort.
\item
The class of regular simple graphs is not $\CMSOL$-definable.   For a general method for
proving non-definability in $\CMSOL$, see \cite{makowsky2014connection}. Hence Theorem \ref{th:BS-CMSOL} cannot be applied
to the sequence $A295193$, which counts the  number of regular simple graphs on $n$ labelled nodes. 
In contrast to this, $r$-regular graphs are $\FOL$-definable, hence Theorem \ref{th:BS-CMSOL} can be applied
easily to the sequence $RG(n,r)$ which counts the number of labelled $r$-regular graphs.
The existence of recurrences for fixed $r$ is discussed in \cite{mckay1983applications} and the references cited therein.
For $r=2,3$ this is  $A110040$.
Recurrences for $r = 0, 1, 2$ are found easily. For $r=3,4$ explicit recurrences were published
in \cite{read1970some,read1980number}, and for $r=5$ in \cite{goulden1983hammond}. 
The recurrence for $r=5$ is linear but very long.
In \cite{gessel1990symmetric}, it is shown that $RG(n,r)$ is holonomic (P-recursive) for every $k \in \N^+$.
We have not checked whether $RG(n,r)$ is holonomic over $\Z$.
In \cite{read1980number} it is shown that $RG(n,4)$ is polynomially recursive, but the equations given there
do not show that $RG(n,4)$ is polynomially recursive over $\Z$.
%It seems likely but still unproven that for general fixed $r$ there are also linear recurrences for $RG(n,r)$, hence that
%$RG(n,r)$ is holonomic. 
It seems that 
Theorem \ref{th:BS-CMSOL} is the most suitable method to show that for each $r$ the sequence $RG(n,r)$ is MC-finite.
\end{enumerate}
\end{examples}

We will use an extension to $\CMSOL$, $\MSOL$ extended by modular counting quantifiers,
from \cite{fischer2011application}, and a new extension which allows the use of hard-wired constants 
and is described in Section \ref{se:constants}.

Clearly, $S_{\phi}(n)$ is computable by brute force, given $\phi$ and $n$.
In \cite{specker1990application}, it is mentioned that
$S^m_{\phi}(n) = S_{\phi}(n) \bmod{m}$ can be computed more efficiently, but no details are given.
Only the special case for $Q^m(n)$ is given, where $Q(n)$ is the number of quasi-orders on $[n]$.
%\begin{problem}
%\label{problem-1}
%Given $\phi$ and $m$ find algorithms for computing $S_{\phi}(n)$ and  $S^m_{\phi}(n)$  and
%determine upper and lower bounds for them.
%\end{problem}

%-------------------------------------------------------------------
%\input{rj-immediate}
\section{Immediate consequences of the Specker-Blatter Theorem}
\label{se:immediate}
\subsection{The Bell numbers $B(n)$}

The Bell numbers $B(n)$ count the number of partitions of the set $[n]$.
This is the same as counting the number of equivalence relations on $[n]$,
which is expressible by an $\FOL$-formula. 
Therefore, we get immediately from Theorem \ref{th:BS-CMSOL} that:
\begin{theorem}
\label{th:bell}
The Bell numbers $B(n)$ are MC-finite.
\end{theorem}
The Bell numbers do satisfy some known congruences.
For $m=p$ a prime, they satisfy the Touchard congruence 
$$
B(p+n) \equiv B(n) + B(n+1) \mod{p}.
$$
However, this is not enough to establish MC-finiteness.

The Bell numbers are not C-finite, because they grow too fast. 
The following estimate is due to \cite{de1981asymptotic,berend2010improved}.
\begin{proposition}
\label{prop:b-growth}
For every $n \in \N^+$
$$ \left(\frac{n}{e \ln n}\right)^n \leq B(n).$$ 
Furthermore, for every $\epsilon > 0$ there is $n_0(\epsilon)$ such that for all $n \geq n_0(\epsilon)$ 
$$B(n) \leq \left(\frac{n}{e^{1-\epsilon} \ln n}\right)^n.$$
\end{proposition}
Better estimates are known, see \cite[Proposition VIII.3]{flajolet2009analytic}, but are not needed here.
Another way to see that Bell numbers are not C-finite is by noticing that they are not holonomic,
\cite{klazar2003bell}. There, and in \cite{banderier2002generating}, some variations of Bell numbers are also studied:
\begin{definition}
\begin{enumerate}[(i)]
\item
$B(n)_{k,m}$ counts the number of partitions of $[n]$ which have $k$ blocks modulo $m$.
\item
$B(n)^{\pm} = B(n)_{0,2} - B(n)_{1,2}$ which  are the Uppuluri-Carpenter numbers $A000587$.
\item
$B(n)^{bc}$ counts the number of {\em bicolored partitions} of $[n]$, i.e., the partitions of $[n]$
where the blocks are colored with two non-interchangeable colors $C_1, C_2$,
$A001861$.
\end{enumerate}
\end{definition}
\begin{theorem}
\label{th:klazar}
The sequences $B(n), B(n)_{k,m}, B(n)^{\pm}, B(n)^{bc} $ are not holonomic, hence not C-finite, but they are
MC-finite.
\end{theorem}
\begin{proof}
That they are not holonomic is shown in \cite{klazar2003bell}, and in \cite{banderier2002generating}.
To see that they are MC-finite, we apply Theorem \ref{th:BS-CMSOL}.
\begin{enumerate}[(i)]
\item
$B(n)_{k,m}$ is definable in $\CMSOL$. We say that there is a set $X \subseteq [n]$ which intersects every block in
exactly one element, and $|X| = k \mod{m}$.
\item
$B(n)^{\pm}$ is the difference of two MC-finite sequences, hence MC-finite.
\item
$B(n)^{bc}$ counts the number of binary and unary relations $E, C_1, C_2$ on $[n]$ such that
$E$ is an equivalence relations, $C_1, C_2 \subseteq [n]$ partition $[n]$, and each of them
is closed under $E$.
\end{enumerate}
\end{proof}

\subsection{Counting transitive relations}
The Bell numbers $B(n)$ count the number of equivalence relations $E(n)$ on a set $[n]$.
Similarly we can look at 
the number of linear quasi-orders (linear pre-orders) $LQ(n)$,
the number of quasi-orders (pre-orders) $Q(n)$,
the number of partial orders $P(n)$, and 
the number of transitive relations $T(n)$ on the set $[n]$.
%The numbers $LQ(n)$ are called {\em ordered Bell numbers} or {\em Fubini numbers}, often denoted 
%in the literature by $a(n)$.
These integer sequences were analyzed in \cite{pfeiffer2004counting}.
They are all definable in $\FOL$, and we have
\begin{proposition}
\label{transitive}
$
B(n)=E(n) \leq LQ(n) \leq P(n) \leq Q(n) \leq T(n).
$
\end{proposition}
\begin{proof}
$E(n) \leq LQ(n)$: We can turn an equivalence relation into a linear quasi-order
by linearly ordering the equivalence classes.

$LQ(n) \leq P(n)$: Each linear quasi-order can be made into a partial order 
by replacing every set of mutually equi-comparable elements in a linear quasi-order with an anti-chain.

$P(n) \leq Q(n)$: Each partial order is also a quasi-order.
%Each quasi-order can be made into a partial order by linearly ordering
%the sets of elements which are mutually bi-comparabel.

$Q(n) \leq T(n)$: Each quasi-order is transitive.
\end{proof}
Hence we get using the Specker-Blatter Theorem and Proposition \ref{transitive}:
\begin{theorem}
\label{th:1}
The sequences 
$
B(n)=E(n), LQ(n),  P(n),  Q(n)$  and 
$T(n) $
are MC-finite but not C-finite.
\end{theorem}

\subsection{Stirling numbers of the second kind}
Let $S(n,k)$ be the number of partitions of $[n]$ into $k$ non-empty blocks.
$S(n,k)$ is also known as the Stirling number of the second kind. Clearly,
$$
B(n) = \sum_k S(n,k).
$$
\begin{theorem}
\label{th:stirling}
For fixed $k=k_0$ the sequence $S(n, k_0)$ is C-finite, and hence MC-finite.
\end{theorem}
This can be seen by observing that 
$S(n, k_0)$ has a rational generating function, see \cite[7.47]{graham1989concrete}.
%\footnote{
%thanks to Ira Gessel for pointing this out.
%}
$$
\sum_{n=0}^\infty S(n,k_0) x^n = \frac{x^{k_0}}{(1-x)(1-2x)\cdots (1-k_0x)}.
$$

\subsection{Lah numbers $Lah(n)$, $A001286$}
%AAA
%The Stirling numbers of the second kind $S(n,k)$ give the number of partitions of a set with $n$ elements into 
%$k$ nonempty subsets.
If we modify the Stirling numbers of the second kind $S(n,k)$ such that the elements in the blocks of the 
partition are ordered between them, 
we arrive at the somewhat less known Lah number $Lah(n,k)$, $A001286$, 
introduced by I. Lah in \cite{lah1954new,lah1955neue} in the context
of actuarial science. 
Good references for Lah numbers are \cite{graham1989concrete,charalambides2018enumerative}.
The Lah numbers are also coefficients expressing rising factorials $x^{(n)}$ in terms of 
falling factorials $x_{(n)}$.
\begin{proposition}
\label{pr:Lah-identity}
$$
x^{(n)} = \sum_{k=1} Lah(n,k) x_{(k)}
\text{  and  }
x_{(n)} = \sum_{k=1} Lah(n,k) x^{(k)}
$$
\end{proposition}
In \cite{guo2015six} six proofs of Proposition \ref{pr:Lah-identity} are given.
Furthermore, $Lah(n) = \sum_k Lah(n,k)$.

$Lah(n)$ counts the number of linear quasi-orders on $[n]$, hence $Lah(n) = LQ(n)$, and
$Lah(n,k)$ counts the number of linear quasi-orders on $[n]$ with $k$ sets of equi-comparable elements.
Two elements $u,v$ in a quasi-order are {\em equi-comparable} if both $u \leq v$ and $v \leq u$.
This is again definable in first order logic $\FOL$.

There are explicit formulas:
\begin{proposition}
\label{prop:Lah}
\begin{gather}
Lah(n,k) = \frac{n!}{k!} \cdot {{n-1} \choose {k-1}} =  \sum_{j=0}^n s(n,j) S(j,k) 
\label{eq:lah-1}
\\
Lah(n) = \sum_k Lah(n,k) = n! \sum_k \frac{1}{k!} \cdot {{n-1} \choose {k-1}} 
\label{eq:lah-2}
\end{gather}
where $s(n,j)$ are the Stirling numbers of the first kind, see \cite{comtet2012advanced}.
\end{proposition}
There is also a recurrence relation:
\begin{gather}
Lah(n+1,k) = Lah(n, k-1) + (n+k) Lah(n,k)
\label{eq:lah-3}
\end{gather}
But again this is not enough to establish C-finiteness or MC-finiteness, since it is a recurrence
involving both $n$ and $k$.

\begin{theorem}
\label{th:lah}
Both $Lah(n)$ and $Lah(n,k_0)$ are MC-finite but not C-finite.
\end{theorem}
\begin{proof}
It follows directly from Equation (\ref{eq:lah-1}), and also from Equation (\ref{eq:lah-3}), 
%and (\ref{eq:lah-2}) 
that for $k=k_0$ fixed the sequence $Lah(n,k_0)$ is not
C-finite.
%----------------------------
MC-finiteness again follows using Theorem 
\ref{th:BS-CMSOL}.
%\ref{th:SP-fol}.
\end{proof}
Note however that the recurrence relation given in Equation (\ref{eq:lah-3})
does not have constant coefficients.

\subsection{Summary so far}
Table \ref{table-1} summarizes the results which are direct consequences of the growth arguments or non-holonomicity (NH)
and the Specker-Blatter Theorem \ref{th:BS-CMSOL} (SB).

\begin{center}
\begin{table}[h]
\begin{tabular}{||l|l|l|l||l|l|l|}
\hline
\hline
Series & C-finite & Proof & Theorem & MC-finite & Proof & Theorem \\
\hline
\hline
$S(n) = B(n)$ & no & Growth & \ref{th:1} &  yes & SB  & \ref{th:bell} \\
\hline
$S(n,k_o)$ & yes & gen.fun& \ref{th:stirling} &    yes & gen.fun &\ref{th:stirling} \\ %or \ref{th:c-fin} \\
\hline
$B(n)^{\pm}$ & no & NH & \ref{th:klazar} &    yes & SB &\ref{th:klazar} \\ %or \ref{th:c-fin} \\
\hline
$B(n)^{bc}$ & no & NH & \ref{th:klazar} &    yes & SB &\ref{th:klazar} \\ %or \ref{th:c-fin} \\
\hline
\hline
$LQ(n)$ & no & Growth & \ref{th:1} &   yes & SB &  \ref{th:1} \\
\hline
$Q(n)$ & no & Growth & \ref{th:1} &   yes & SB &  \ref{th:1} \\
\hline
$P(n)$ & no & Growth & \ref{th:1} &   yes & SB &  \ref{th:1} \\
\hline
$T(n)$ & no & Growth & \ref{th:1} &   yes & SB &  \ref{th:1} \\
\hline
\hline
$Lah(n) = LQ(n)$ & no & Growth & \ref{th:lah} & yes & SB & \ref{th:lah} \\
\hline
$Lah(n, k_0)$ & no & Growth & \ref{th:lah} &   yes & SB & \ref{th:lah} \\
\hline
\hline
\end{tabular}
\caption{Direct consequences of the Specker-Blatter Theorem}
\label{table-1}
\end{table}
\end{center}
%-------------------------------------------------------------------
%\input{rj-restricted}
\section{Restricted set partitions}
\label{se:restricted}
%\section{Restricted Stirling and Lah numbers}
The new results of this paper concern C-finiteness and MC-finiteness
for restricted versions of set partitions.
%for restricted versions of Stirling and Lah numbers.
We have two kinds of restrictions in mind.
The first are {\em positional restrictions} which impose conditions on the positions of the elements of $[n]$ where $[n]$
is equipped with its natural order.
The second are {\em size restrictions} which impose conditions on the size of the blocks or the number of the blocks.

\subsection{Global positional restrictions}

\begin{definition}
Let $A$ and $B$ be two blocks of a partition of $[n]$.
\begin{enumerate}[(i)]
\item
$A$ and $B$ are {\em crossing} if there are elements
$a_1, a_2 \in A$ and $b_1, b_2 \in B$ such that 
$a_1 < b_1 < a_2 < b_2$ or
$b_1 < a_1 < b_2 < a_2$.
\item
Let 
$\min{A},\max{A},\min{B},\max{B}$ 
the smallest and the largest elements in $A$ and $B$.
$A$ and $B$ are {\em overlapping} if 
$\min{A} < \min{B} < \max{A} < \max{B}$ or
$\min{B} < \min{A} < \max{B} < \max{A}$.
\item
If $A$ and $B$ are overlapping they are also crossing, but not conversely.
\item 
The number $B(n)^{nc}$ of 
non-crossing set partitions on $[n]$ is one of the interpretations of the Catalan numbers, \cite{roman2015introduction}.
\item
The Bessel number
$B(n)^B$ ($A006789$) is the number 
of non-overlapping set partitions on $[n]$, \cite{flajolet1990non}.
\end{enumerate}
\end{definition}
The Catalan numbers $C(n)$ are not holonomic and not MC-finite.
In \cite{banderier2002generating} it is shown that the Bessel numbers $B(n)^B$ are not holonomic.
%\begin{problem} \label{prob:bessel} 
Are the Bessel numbers $B(n)^B$ MC-finite?  
%\end{problem}
The positional restrictions here are {\em global} in the sense that they involve all of the elements
of $[n]$ with their natural order.
For non-holonomic integer sequences $s(n)$ that count the number of set partitions subject to global positional restrictions,
we have currently no tools to decide whether they are MC-finite or not.

Next, we look at {\em local} positional restrictions 
one can impose on
Stirling and Lah numbers, \cite{broder1984r,wagner1996generalized,nyul2015r,benyi2019some,benyi2019restricted}.
They are local because they only put restrictions on the positions of a fixed number
of elements of $[n]$ with their natural order.

\subsection{Local positional and size restrictions}
Recall that we denote by $[n]$ the set $\{1, 2. \ldots, n\}$.
We denote by $S_r(n,k)$ the number of partitions of $[n+r]$ into $k+r$ non-empty blocks
with the additional condition that the first $r$ elements are in distinct blocks.
The elements $1, \ldots, r$ are called {\em special elements} and the partitions where the first $r$ elements
are in distinct blocks are called $r$-partitions. When dealing with definability we view the special elements
as {\em hard-wired} constants, i.e., constant symbols $a_i, 1 \leq i \leq r$ 
with a fixed interpretation by elements of $[n+r]$.

We define $S_r(n)= B_r(n)$ by
$$
S_r(n) = \sum_k S_r(n,k).
$$
$Lah_r(n,k)$, $A143497$, and $Lah_r(n)$ are defined analogously, 
%with the additional condition that $a_1 < a_2 < \ldots < a_r$, 
with the condition that $a_1 < a_2 < \ldots < a_r$ are in different blocks.
\cite{nyul2015r,shattuck2016generalized}.

Let $A \subseteq \N$. We denote by
$S_{A,r}(n) = B_{A,r}(n)$, 
$S_{A,r}(n,k)$,
$Lah_{A,r}(n)$ and $Lah_{A,r}(n,k)$
the number of corresponding partitions where every block has its size in $A$. 

For $r=0$, in the absence of special elements, we just write $S_A(n) = B_A(n)$, $S_A(n,k)$,
$Lah_A(n)$ and $Lah_A(n,k)$.

A set $A \subseteq \N$ is {\em (ultimately) periodic} if there exist $p,n_0 \in \N^+$ such that for all $n \in \N$
($n \geq n_0$)
we have $n \in A$ iff $n+p \in A$.
In other words, the characteristic function $\chi_A(n)$ of $A$ is ultimately periodic in the usual sense,
$\chi_A(n) = \chi_A(n+p)$ ($n \geq n_0$).
Analogous definitions can be made for $LQ(n)$, denoted by $LQ_{A,r}$, and also called $r$-Fubini sequences, 
with OEIS-number $A232472$.

\subsection{Main results for restricted set partitions}
Our  results for restricted set partitions are summarized in Tables
\ref{table-2}, 
\ref{table-3},
\ref{table-4} and 
\ref{table-5} 
below.
%\\
FM refers to the proof method of \cite{fischer2008linear,fischer2011application}.
%\\
%SB+ refers to the extension of the Specker-Blatter Theorem to definability in $\CMSOL$,
SB* refers to the extension of the Specker-Blatter Theorem to allow a fixed finite set of special elements as
hard-wired constants.

%\newpage
%\input{rj-t2}
\begin{center}
\begin{table}[h]
\begin{tabular}{|l||l|l|l||l|l|l|}
\hline
\hline
Series & C-finite & Proof & Theorem & MC-finite & Proof & Theorem \\
\hline
\hline
$S_A(n) = B_A(n)$ & no & Growth & \ref{th:g-bell} & yes & SB* & \ref{th:A-periodic} \\
\hline
$S_A(n,k_0)$ & yes & FM & \ref{th:c-fin} & yes &  FM &  \ref{th:c-fin} \\
\hline
$Lah_A(n) = LQ_A(n)$ & no & Growth & \ref{th:g-lah} & yes & SB* & \ref{th:A-periodic} \\
\hline
$Lah_A(n, k_0)$ & no & Growth & \ref{th:g-lah} & yes & SB* & \ref{th:A-periodic} \\
\hline
\hline
\end{tabular}
\caption{With ultimately periodic $A$ only}
\label{table-2}
\end{table}
\end{center}
\begin{center}
\begin{table}[h]
\begin{tabular}{|l||l|l|l||l|l|l|}
\hline
\hline
Series & C-finite & Proof & Theorem & MC-finite & Proof & Theorem \\
\hline
\hline
$S_r(n)=B_r(n)$ & no & Growth & \ref{th:g-bell} & yes & SB* &  \ref{co:SB*} \\
\hline
$S_r(n, k_0)$ & yes & FM &\ref{th:c-fin} &  yes & FM &  \ref{th:c-fin} \\
\hline
\hline
$Lah_r(n, k_0)$ & no & Growth &\ref{th:g-lah} &  yes & SB* & \ref{co:SB*} \\
\hline
\hline
\end{tabular}
\caption{With hard-wired constants only}
\label{table-3}
\end{table}
\end{center}
\begin{center}
\begin{table}[h]
\begin{tabular}{|l||l|l|l||l|l|l|}
\hline
\hline
Series & C-finite & Proof & Theorem & MC-finite & Proof & Theorem \\
\hline
\hline
$S_{A,r}(n)= B_{A,r}(n)$ & no & Growth & \ref{th:bell-A} &  yes & SB* &   \ref{co:SB*} \\
\hline
$S_{A,r}(n, k_0)$ & yes & FM & \ref{th:c-fin}  & yes & FM &  \ref{th:c-fin} \\
\hline
$Lah_{A,r}(n, k_0)$ & no & Growth & \ref{th:bell-A} & yes & SB* & \ref{co:SB*} \\
\hline
\hline
\end{tabular}
\caption{With ultimately periodic $A$ and hard-wired constants }
\label{table-4}
\end{table}
\end{center}
These results also hold for $LQ_{A,r}$, the $r$-Fubini numbers, and other similarly defined sequences.

\begin{center}
\begin{table}[h]
\begin{tabular}{|l||l|l|l||l|l|l|}
\hline
\hline
Series & C-finite & Proof & Theorem & MC-finite & Proof & Theorem \\
\hline
\hline
$B(n)^B$ & no & NH & \cite{banderier2002generating} &  ??? & ??? &  ---  \\
\hline
$B(n)^{nc} = C(n)$ & no & NH & \cite{roman2015introduction} &  no & \cite[Theorem 9.4]{roman2015introduction} 
 & \cite{banderier2002generating} \\

\hline
%-----
\hline
\hline
\end{tabular}
\caption{With global positional restrictions}
\label{table-5}
\end{table}
\end{center}
%-------------------------------------------------------------------
%\input{rj-rproofs}
\section{Proofs for the restricted cases}
\label{se:rproofs}
For the analysis of MC-finiteness in the restricted cases we need some additional tools.

\subsection{Ultimate periodicity of $A$}
\label{se:periodic}
\ifmargin
\marginpar{s-periodic}
\else\fi

Recall that a formula with a {\em modular counting quantifier $C_{b,m}x \phi(x)$} is true in a structure $\fB$ if
the cardinality of the set of elements in $\fB$ which satisfy $\phi(x)$, satisfies 
$$|\{a \in B : \phi(a) \}| \equiv b \mod{m}.$$
$\CMSOL$ is the logic obtained from $\MSOL$ by extending it with all the modular counting
quantifiers $C_{b,m}$.
In \cite{pr:FischerMakowsky03} the Specker-Blatter Theorem was extended to hold for $\CMSOL$, as already
stated in Theorem \ref{th:BS-CMSOL}.
$\CMSOL$ is also needed to prove the following lemma:

\begin{lemma}
\label{le:periodic}
Let $A$ be ultimately periodic and $\psi(x)$ be a formula of $\CMSOL$.
Then there is a sentence $\psi_A \in \CMSOL$ such that
in every finite structure $\fB$ we have
$$\fB \models \psi_A \text{ iff }|\{b \in B : \psi(b) \}| \in A$$
\end{lemma}
%--------------------------------------------------------
\begin{proof}
If  $A = A_{a,m}= \{ n \in \N: n \equiv a \mod m \}$ the formula $\psi_A$
is the sentence $C_{a,m} x \psi(x)$.

Next we observe that if $A$ is ultimately periodic there are finitely many $a_1, \ldots, a_k$ and $q$
such that $A = \bigcup_{i=0}^k A_i$ with $A_0 \subseteq [q]$ and 
$A_i = \{ n > q: n \equiv a_i \mod m \}$.
We proceed in steps:
%----------------Eldar's version--------------
\begin{enumerate}[(i)]
\item
$\exists^{\geq k}x\psi(x):=\exists x_1,\ldots,x_k \bigwedge_{i=1}^k\psi(x_i)\wedge\bigwedge_{1\leq i<j\leq k}(x_i\neq x_j)$
says that there are at least $k$ elements that satisfy $\psi(x)$.
\item
$\exists^{=k}x\psi(x):=\exists^{\geq k}x\psi(x)\wedge\neg\exists^{\geq k+1}x\psi(x)$
says that there are exactly $k$ such elements.
\item
$\psi{A_0}:=
(\exists^{<q} x \psi(x)
\rightarrow
\bigvee_{j\in A_0}\exists^{=j}x\psi(x))$
says that if the number of elements satisfying $\psi(x)$ is less than $q$
then the number of such elements has exactly one of the cardinalities in $A_0$.
\item
$\psi{A_i}:=\exists^{\geq q+1}x\psi(x)\wedge C_{a,m}x\psi(x)$
says that if the number of elements satisfying $\psi(x)$ is bigger or equal than $q$
then the number such elements equals $a_i \bmod{m}$.
\item
$\psi(x)_A := \bigvee_{i=0}^k \psi{A_i}(x)$ is the required formula.
\end{enumerate}
%-----------------end Eldar's version------------------
\end{proof}

Theorem \ref{th:BS-CMSOL} together with Lemma \ref{le:periodic} gives immediately:
\begin{theorem}
\label{th:A-periodic}
Assume that $A$ is ultimately periodic.
Then the sequences
$B_A(n)= S_A(n), Lah_A(n)$ and $Lah_A(n, k_0)$ are MC-finite.
\end{theorem}

\subsection{Growth arguments}

We first discuss growth arguments for $B_A(n)= S_A(n), Lah_A(n)$ and $Lah_A(n, k_0)$.
\begin{theorem}
\label{th:bell-A}
Let $A \subseteq \N$ be infinite and ultimately periodic. Then 
$B_A(n)= S_A(n), Lah_A(n)$ and $Lah_A(n, k_0)$ are not C-finite.
\end{theorem}
\begin{proof}
First we prove it for $B_A(n)$ and $A = A_m = \{ n \in \N: n \equiv 0 \mod m \}$. 
Let $P_1, \ldots , P_k$ be a partition of $[n]$.
We replace in each $P_i$ every element by $m$ elements. This gives us a partition of $[mn]$ with each block
of size in $A_m$. Hence 
$$
P_A(mn) \geq P(n) \geq \left(\frac{n}{e \ln n}\right)^n
$$
or, equivalently,
$$
P_A(n) \geq P(n/m) \geq \left(\lfloor\frac{n/m}\rfloor{e \ln \lfloor n/m \rfloor}\right)^{\lfloor n/m \rfloor}
$$
which still grows superexponentially.

Next we assume that $A = A_{k,a,m}= \{ n \in \N: n \equiv a \mod m, n \geq k \}$.
We proceed as before, but additionally add $mr+a$ elements to each block, for $r$ large enough.
%---------------------
%only instead of adding a you need to add mr+a for r large enough that this goes into the periodic part of A.
%-------------------------------------------------
Finally, we note that for every infinite (ultimately) periodic set $A$  
there is a set $A_{k, a,m}$ for some $k, a, m \in \N^+$ such that 
$ A_{k,a,m} \subseteq A$
%there is a set $A_{a,m}$ for some $a, m \in \N^+$ such that 
%$ A \cap A_{a,m} =B$ is infinite. 
%We then repeat the argument with $B$ instead of $A$.

For $Lah_A(n)$ and $Lah_A(n, k_0)$ we proceed similarly using Proposition \ref{prop:Lah}.
\end{proof}

Next we discuss growth for 
$Lah(n,k_0)$, $Lah(n) = \sum_k Lah(n,k)$ and $Lah_r(n, k_0)$.

We have seen in Proposition \ref{prop:b-growth} that 
$$
\left(\frac{n}{e \ln n}\right)^n \leq B(n) \leq \left(\frac{n}{e^{1-o(1)} \ln n}\right)^n
$$

We now show
\begin{lemma}
\label{le:bell-3}
$B_r(n) \geq B(n)$
\end{lemma}
\begin{proof}
Every partition of $[n]$ gives rise to at least one partition of $[n+r]$ where the first $r$ elements are in
distinct blocks containing only one element.
\end{proof}

From  Proposition \ref{prop:c-finite}, \ref{transitive} 
and \ref{le:bell-3}
%\ref{le:bell-1} 
%\ref{le:bell-2} 
we get:

\begin{theorem}
\label{th:g-bell}
The sequences $B(n)$ and $B_r(n)$
%LQ(n), Q(n), P(n)$ and $T(n)$ 
are not C-finite.
\end{theorem}

\begin{lemma}
\label{le:lah-1}
For $k_0, r$ fixed, the Lah number $Lah(n,k_0)$ satisfy  the following:
\begin{enumerate}[(i)]
\item
$Lah(n,k_0) = {{n-1} \choose {k_0-1}} \frac{n!}{k_0!}$,
\item
$Lah(n) \geq Lah(n, k_0)$, and
\item
$Lah_r(n, k_0) \geq Lah(n, k_0)$.
\end{enumerate}
\end{lemma}
\begin{proof}
(i) is from \cite{lah1954new,lah1955neue}. (ii) follows from (i), and (iii) is proved like Lemma \ref{le:bell-3}.
\end{proof}

This gives immediately

\begin{theorem}
\label{th:g-lah}
Let $k_0$ be fixed. The sequences $Lah(n,k_0)$, $Lah(n) = \sum_k Lah(n,k)$ and $Lah_r(n, k_0)$ are not C-finite.
\end{theorem}

\subsection{Hard-wired constants}
%\label{se:hard}
%\ifmargin
%\marginpar{s-hard}
%\else\fi

Recall that a constant is {\em hard-wired} on $[n]$ if its interpretation is fixed.

The Specker-Blatter  Theorem is originally proved for classes of structures with
a finite number of binary relations. 
It is false for one quaternary relations \cite{ar:Fischer02}.
It was announced recently that it is also false for one ternary relation, \cite{ar:FischerMakowsky2022}.

The Specker-Blatter Theorem remains true when adding a finite number of unary relations.
This is so because a unary relation $U(x)$ can be expressed as a binary relation
$R(x,x)$ which is false for $R(x,y)$ when $x \neq y$.

Adding constants comes in two flavors, with variable interpretations, or hard-wired.
Assume we want to count the number of unary predicates $P$ on $[n]$ which contain the interpretation 
of a constant symbol $c$.
There are $n$ possible interpretations for $c$ and $2^{n-1}$ interpretations for sets not containing $c$,
hence $n2^{n-1}$ many such sets.
However, if $c$ is hard-wired to be interpreted as $1 \in [n]$, there are only $2^{n-1}$ many such sets.

Constants can be represented as unary predicates the interpretation of which is a singleton.
If we do this, the Specker-Blatter Theorem holds, but we cannot model the $r$-Bell numbers like this.
To prove that the $r$-Bell numbers are MC-finite one has to deal with $r$ many hard-wired constants.
Adding a finite number of hard-wired constants needs some work.
In Appendix \ref{se:constants} we show how to eliminate a finite number of hard-wired constants
for the case of $S_r(n)$. The proof generalizes. In \cite{ar:FischerMakowsky2022} the more general
version is proved:

\begin{theorem}
\label{th:SB*}
Let $\tau_r$ be a vocabulary with finitely many binary and unary relation symbols, and $r$ hard-wired constants.
Let $\phi$ be a formula of 
$\CMSOL(\tau_r)$.
Then $S_{\phi}(n)$ is MC-finite.
\end{theorem}
\begin{corollary}
\label{co:SB*}
The sequences
$S_r(n)=B_r(n)$, $Lah_r(n, k_0)$, $S_{A,r}(n)= B_{A,r}(n)$,  $Lah_{A,r}(n, k_0)$
are MC-finite.
\end{corollary}

%-------------------------------------------------------------------
%\input{rj-fm}
\subsection{Proving C-finiteness}
\label{se:fm}
\ifmargin
\marginpar{rj-fm}
\else\fi

In this subsection we explain a special case of the method used in \cite{fischer2008linear}
to prove C-finiteness.
It is based on counting partitions of graphs satisfying additional properties and computing these partitions for 
iteratively constructed graphs.

\subsubsection{Counting partitions with a fixed number of blocks}
Let $G = (V(G), E(G))$ be a graph, and $k_0 \in \N$.  We look at partitions $P_1(G), \ldots, P_{k_0}(G)$ 
of $V(G)$ which can be described
in first order logic $\FOL$.  The following are three typical examples:
\begin{examples}
\begin{enumerate}[(i)] 
\item
The underlying sets of $G[P_i(G)]$ form a partition of $V(G)$ without further restrictions.
\item
For each $i \leq k_0$ the induced graph $G[P_i(G)]$ is edgeless (proper coloring).
\item
Let $\mathcal{P}$ be a graph property.
For each $i \leq k_0$ the set $G[P_i(G)]$  is in $\mathcal{P}$ ($\mathcal{P}$-coloring).
\end{enumerate}
\end{examples}

We look at the counting function
$$
f_{\phi}(G) =
| \{ P_1(G), \ldots, P_{k_0}(G) : \phi(P_1(G), \ldots, P_{k_0}(G)) \}|
$$
defined using an $\FOL$-formula $\phi$.

Let $A \subseteq \N$ be an ultimately periodic set.
We also look at the  restricted counting function
$$
f_{\phi, A}(G) =
| \{ P_1(G), \ldots, P_{k_0}(G) : \phi(P_1(G), \ldots, P_{k_0}(G)) \text{  and  } |P_i(G)| \in A \}|.
$$

We also allow graphs with a fixed number of distinct 
vertices, which may appear in the formula $\phi$.

\subsubsection{Iteratively constructed graphs}
\begin{definition}
A $k$-colored graph is a graph $G$ together with $k$ sets 
$V_1,V_2,...,V_k\subseteq V(G)$ such that $V_i \cap V_j =\emptyset $ for $i\neq j$.
A basic operation on $k$-colored graphs is one of the following:
\begin{itemize}
\item $Add_i$: add a new vertex of color $i$ to $G$.
\item $Recolor_{i,j}$: recolor all vertices with color $i$ to color $j$ in $G$.
\item $Uncolor_{i}$: remove the color of all vertices with color $i$. Uncolored vertices cannot be recolored again.
\item $AddEdges_{i,j}$: add an edge between every vertex with color $i$ and every vertex with color $j$ in $G$.
\item $DeleteEdges_{i,j}$: delete all edges between vertices with color $i$ and vertices with color $j$ from $G$.
%\item $Duplicate$: Add a disjoint copy of $G$ to $G$.
\end{itemize}
A unary operation $F$ on graphs is elementary if $F$ is a finite composition of basic 
operations on $k$-colored graphs (with $k$ fixed). 
We say that a sequence of graphs $\{G_n\}$ is iteratively constructed 
if it can be defined by fixing a graph $G_0$ and defining $G_{n+1}=F(G_n)$ for an elementary operation $F$.
\end{definition}

\begin{example}
\label{exIterativelyConstructed}
The following sequences are iteratively constructed:
\begin{itemize}
\item The complete graphs $K_n$ can be constructed using two colors: 
Fix $G_0$ to be the empty graph, and the operation $F$, given a graph $G_n$, adds a vertex with color 2, 
adds edges between all vertices with color 2 and color 1, and recolors all vertices with color 2 to color 1.
\item The paths $P_n$ can be constructed using 3 colors: 
Fix $G_0$ to be the empty graph, and the operation $F$, given a graph $G_n$, 
adds a vertex with color 3, adds edges between all vertices with colors 2 and 3, 
recolors all vertices with color 2 to color 1, and recolors all vertices with color 3 to color 2.
\item The cycles $C_n, n \geq 3$ can be constructed by first constructing a path $P_n$ where the first and the last element
have colors $1$ and $2$ different from the remaining vertices. Then we connect the first and last element of $P_n$
by an edge. This needs $5$ colors, but is not iterative.
To make it an iterative construction we proceed as follows.
Given a cycle $C_n$ with with two neighboring vertices of color $1$ and $2$, uncolor all the other vertices and
remove the edge $(1,2)$. Then add a new vertex with color $3$, make edges $(1,3)$ and $(3,2)$, uncolor the 
old vertices colored by $1$, 
and then recolor $3$ to have color $1$.
\end{itemize}
\end{example}

\begin{remark}
In \cite{fischer2008linear} 
there was an additional operation allowed
\begin{itemize}
\item $Duplicate$: Add a disjoint copy of $G$ to $G$,
\end{itemize}
assuming erroneously that $Duplicate$ behaves like a unary operation on graphs.
Although it looks like a unary operation on graphs, the sequence of graphs
$$
G_0 = E_1, G_{n+1} = Duplicate(G_n)
$$
grows too fast and does not fit the framework that the authors have envisaged in \cite{fischer2008linear}.
\end{remark}

\subsubsection{The FM method}

In this framework \cite{fischer2008linear} proved the following:

\begin{theorem}[The Fischer-Makowsky Theorem]
\label{th:FM}
Let $G_n$ be an iteratively constructed sequence of graphs, $A \subseteq \N$ be ultimately periodic,  and 
$$
f_{\phi}(G_n) =
| \{ P_1(G_n), \ldots, P_{k_0}(G_n) : \phi(P_1(G_n), \ldots, P_{k_0}(G_n)) \}|
$$
and
$$
f_{\phi, A}(G_n) =
| \{ P_1(G_n), \ldots, P_{k_0}(G_n) : \phi(P_1(G_n), \ldots, P_{k_0}(G_n)) \text{  and  } |P_i(G_n)| \in A \}|,
$$
where $\phi \in \CMSOL$.
Then the sequences $f_{\phi}(G_n)$ 
and
$f_{\phi, A}(G_n)$ are C-finite.
\end{theorem}

We now use Theorem \ref{th:FM} to prove:

\begin{theorem}
\label{th:c-finite}
\label{th:c-fin}
Let $A$ be ultimately periodic, $r, k_0 \in \N$. Then
$S(n,k_0)$,
$S_A(n,k_0)$,
$S_r(n,k_0)$ and
$S_{A,r}(n,k_0)$
are C-finite.
\end{theorem}
\begin{proof}
It suffices to prove it for $S_{A,r}(n,k_0)$. 
The other cases can be obtained by setting 
$r=0$ and/or $A= \N$.

We have to show that
$S_{A,r}(n,k_0)$
is of the form 
$f_{\phi, A}(G_n)$.

We define an iteratively constructed sequence of graphs $G =(V(G), E(G), v_1, \ldots, v_r)$
with $r$ distinct vertices as follows.
$G_0 = (K_r, v_1, \ldots, v_r)$.
$G_{n+1} = G_n \sqcup K_1$.

Now take $\phi(P_1, \ldots , P_{k_0}, v_1, \ldots , v_r)$
which says that the $P_i$'s form a partition and for each $i \leq r$  
the distinguished vertex $v_i$ belongs to $P_i(G)$.
\end{proof}

Further details are given in Appendix \ref{se:c-finite}.
%-------------------------------------------------------------------
%\input{rj-conclu}
\section{Conclusions and further research}
\label{se:conclu}

In the first part of the paper we introduced MC-finiteness as a worthwhile topic in the study
of integer sequences. We surveyed two methods of establishing MC-finiteness of such sequences.
%Besides proving it adhoc for a specific sequence, we presented two methods.
In Theorem \ref{th:MC}, MC-finiteness follows from the existence of polynomial recurrence relations
with coefficients in $\Z$. In Theorem \ref{th:BS-CMSOL}, MC-finiteness follows from a logical definability
assumption in Monadic Second Order Logic augmented with modular counting quantifiers $\CMSOL$.
We have compared the advantages and disadvantages of the methods, and we have used
the logic method of Theorem \ref{th:BS-CMSOL} to give quick and transparent proofs of MC-finiteness.

%-----------------------------------------
\ifskip\else
\section{Conclusions and suggestions for further work}
\label{se:conclu}
\ifmargin
\marginpar{s-conclu}
\else\fi %margin

The purpose of this paper was to study congruences (MC-finiteness) for
restricted set partition functions.
In the unrestricted case MC-finiteness is a direct consequence of the
Specker-Blatter Theorem from 1981, 
which, unfortunately,
has been widely unnoticed. 
In its original version it only applies to counting labeled structures which use binary relations
definable in Monadic Second Order Logic $\MSOL$.
The direct applications are summarized in Table \ref{table-1}. 
\fi %skip

In the second part of the paper we
got similar results for locally restricted set partition functions
like $B_{A,r}$. For this purpose the Specker-Blatter Theorem has to be extended 
%The first extension extends its applicability to binary relations definable in $\CMSOL$, $\MSOL$
%augmented by modular counting quantifiers. 
in order to 
%The second extension 
count labeled structures where
a fixed number of special elements are in a certain configuration. In the case of $B_{A,r}$,
$A$ is a set of natural numbers and $r$ is a natural number. $B_{A,r}$ counts the number of set partitions
of $[n]$ where the first $r$ elements are in different blocks and
$A$ indicates the possible cardinalities of  the blocks of the partition.
Such an extension is  given
in Theorem \ref{th:SB*}. A proof of a special case of this theorem is given in the appendix.
The general case can be found in \cite{ar:FischerMakowsky2022}.
Our new results are summarized in Tables \ref{table-2}--\ref{table-5}. 

We did not investigate in depth whether 
MC-finiteness 
of the examples 
in Tables \ref{table-2}--\ref{table-5}
can be established directly or by exhibiting suitable polynomial recurrence schemes, in order to apply
Theorem \ref{th:MC}.

\begin{problem}
\label{problem-Bessel}:
Are the Bessel numbers $B(n)^B$ MC-finite?
\end{problem}

\begin{problem}
\label{problem-PRS}:
Find systems of mutual polynomial recurrences for all the examples in Tables \ref{table-2}--\ref{table-4}.
\end{problem}

Instead of set partition functions we can also count the number of, say, partial orders
where
\begin{enumerate}[(i)]
\item
$r$ special elements are in a particular $\CMSOL$ definable configuration,
such as prescribed comparability and incomparability, and
\item
$A$ indicates the possible cardinalities of  certain definable sets, such as antichains or
maximal linearly ordered sets.
\end{enumerate}

Our techniques allow us to show that counting such partial orders on $[n]$ results in MC-finite
sequences.

In \cite{specker1990application} it is suggested that counting the number of quasi-orders $Q^m(n)$
on $[i]$ modulo $m$ is easier than finding the exact value of $Q(n)$. 
%A general version
%of this was stated in the introductions as Problem \ref{problem-1}.

Clearly, $S_{\phi}(n)$ is computable by brute force, given $\phi$ and $n$.
In fact, for $\phi \in \FOL$ the problem is in $\sharp\bP$.
For $\phi \in \CMSOL$ it is in $\sharp\bPH$, the analogue of $\sharp\bP$ for problems definable 
in Second Order Logic, or equivalently,
in the polynomial hierarchy. As noted in \cite[Proposition 11]{makowsky1996arity}, 
there are arbitrarily complex problems in $\bPH$ already definable in $\MSOL$.
However, $S^m_{\phi}(n)$ is in $MOD_m\bP$, respectively in $MOD_m\bPH$, the corresponding modular counting classes
introduced in \cite{beigel1992counting}. 
%It is an open problem, whether $MOD_m\bP$ is strictly weaker that $\sharp\bP$.
It is still open how exactly $MOD_m\bP$ is related to $\sharp\bP$.

In \cite{specker1990application}, it is mentioned that
$S^m_{\phi}(n) = S_{\phi}(n) \mod{m}$ can be computed more efficiently, but no details are given.
Only the special case of $Q^m(n)$ is given, where $Q(n)$ is the number of quasi-orders on $[n]$.

%AAAA Last corrections
\begin{problem}
\label{problem-1}
%\label{problem-bounds}
Given $\phi \in \FOL$ and $m$, find algorithms for computing $S_{\phi}(n)$ and  $S^m_{\phi}(n)$  and
determine upper and lower bounds for them. One may assume that $n$ is encoded in unary.
\end{problem}
\begin{problem}
\label{problem-2}
Same as Problem \ref{problem-1} for $\phi \in \CMSOL$.
\end{problem}
\begin{problem}
Inspired by the remarks above, the following might be a worthwhile project:
\label{problem-3}
Investigate the complexity classes $\sharp\bPH$ and $MOD_m\bPH$ and their mutual relationships.
\end{problem}
%-------------------------------------------------------------------
%\input{rj-oeis-nu}
\section{List of OEIS-sequences}
\label{se:oeis}
\begin{description}
\item[A000108] Catalan numbers $C(n)$.
\item[A000110] Bell numbers $B(n)$.
\item[A000453] Stirling numbers of the send kind $S(n,k)$.
\item[A000587] Uppuluri-Carpenter numbers $A000587$.
\item[A000670] Number of linear quasi-orders (pre-orders) $LQ(n)$.
\item[A000798] Number of quasi-orders (pre-orders) $Q(n)$.
\item[A001035] Number of partial orders $P(n)$.
\item[A001286] Lah numbers $Lah(n)$.
\item[A001861] Bicolored partitions.
\item[A005493] $r$-Bell numbers $B_{A,2}(n)$ for $r=2$.
\item[A005494] $r$-Bell numbers $B_{A,3}(n)$ for $r=3$.
\item[A006905] Number of transitive relations $T(n)$.
\item[A086714] $a(0)=4, a(n+1) = {a(n) \choose 2}$.
\item[A110040] Regular labeled graphs of degree $2$ and $3$.
\item[A143494] $r$-Stirling numbers $S_{A,r}(n,k)$
\item[A143497] $r$-Lah numbers $Lah_{A,r}(n)$.
\item[A232472] $r$-Fubini numbers $LQ_{A,r}$ for $r=2$.
\item[A295193] Regular labeled graphs.
\end{description}

\newpage
%-------------------------------------------------------------------
%\bibliographystyle{plain}
%\bibliography{rj-quotes}
%-------------------------------------------------------------------
%\input{rj-ref}

%-------------------------------------------------------------------
\appendix
\section{More on MC-finiteness}
\label{ap:mc-finite}
\label{se:ProofMC}

\subsection{Polynomial recursive sequences}

A \emph{polynomial recursive sequence}~\cite{cadilhac2021polynomial} is a mutual recurrence in which the recurrence relation is a polynomial. That is, we define $d$ sequences in parallel by initial values $a_1(0),\dots,a_d(0)$ and the recurrence
\[
 a_i(n+1) = P_i(a_1(n),\dots,a_d(n)),
\]
where $P_i$ is a polynomial with rational coefficients. We will only consider recurrences for which $a_i(n) \in \mathbb{N}$ for all $i \in [d]$ and $n \ge 0$.

\begin{theorem}[\cite{cadilhac2021polynomial}] \label{thm:CMPPS20}
Let $m$ be a natural number which is relatively prime to all denominators of coefficients of the defining polynomials $P_1,\dots,P_d$. Then the sequences $a_i(n) \bmod m$ are eventually periodic.
\end{theorem}
\begin{proof}
Notice that
\[
 a_i(n+1) \bmod m = (P_i \bmod m)(a_1(n) \bmod m, \dotsm a_d(n) \bmod m).
\]
Thus the function $P\colon \mathbb{Z}_m^d \to \mathbb{Z}_m^d$ given by
\[
 P(x_1,\dots,x_d) = ((P_1 \bmod m)(x_1,\dots,x_d),\dots, (P_d \bmod m)(x_1,\dots,x_d))
\]
satisfies
\[
 (a_1(n+1) \bmod m,\dots,a_d(n+1) \bmod m) = P(a_1(n) \bmod m,\dots,a_d(n) \bmod m).
\]
Since $\mathbb{Z}_m^d$ is finite, if we start at $(a_1(0) \bmod m, \dots, a_d(0) \bmod m)$ and repeatedly apply $P$, we will eventually enter a cycle.
\end{proof}

This result raises the following question: what happens for other $m$? It turns out that the theorem fails in general for such $m$.

\begin{theorem} \label{thm:main}
Consider the following sequence $A086714$:
%\href{https://oeis.org/A086714}{A086714}:
\[
 a(n+1) = \binom{a(n)}{2}, \quad a(0) = 4.
\]
The sequence $a(n) \bmod 2$ is not eventually periodic.
\end{theorem}

The same result holds (with the same proof) for any $a(0) \ge 4$, 
as well as for any recurrence of the form $a(n+1) = (a(n)+b)(a(n)+c)/2$, 
as long as $b,c$ have different parities and $a(0)$ is chosen so that $a(n) \to \infty$.

\subsection{Proof of Theorem \ref{th:MC}}

Let $\beta(n) = a(n) \bmod 2$. It is not hard to check that the sequence 
$\beta(n) \ldots \beta(n + k - 1)$ depends only on $a(n) \bmod 2^k$. 
It turns out that the opposite holds as well: we can determine $a(n) \bmod 2^k$ from $\beta(n) \ldots \beta(n + k - 1)$.

\begin{lemma} \label{lem:bijection}
Let $a_r,\beta_r$ be defined as above, except with the initial condition $a_r(0) = r$.

For all $k \ge 1$, the function
\[
 \Phi_k(r) = \beta_r(0) \ldots \beta_r(k-1)
\]
is a bijection between $\{0,\dots,2^k-1\}$ and $(0,1)^k$.
\end{lemma}

For example, if $k = 3$, we get the following bijection:
\begin{align*}
\Phi_3(0) &= 000 &
\Phi_3(1) &= 100 &
\Phi_3(2) &= 010 &
\Phi_3(3) &= 111 \\
\Phi_3(4) &= 001 &
\Phi_3(5) &= 101 &
\Phi_3(6) &= 011 &
\Phi_3(7) &= 110
\end{align*}

\begin{proof}
The proof is by induction on $k$. The result is clear when $k = 1$, so suppose $k > 1$.

The first bit of $\Phi_k(r)$ is the parity of $r$, and the remaining bits are $\Phi_{k-1}(s)$, where $s = \binom{r}{2} \bmod 2^{k-1}$. To complete the proof, we show that the mapping $r \mapsto s$ is $2$-to-$1$, with the two pre-images of every $s$ having different parity.

Indeed, suppose that $\binom{a}{2} \equiv \binom{b}{2} \pmod{2^{k-1}}$ for $a,b \in \{0,\dots,2^k-1\}$. Then $a(a-1) \equiv b(b-1) \pmod{2^k}$, and so $2^k \mid a(a-1) - b(b-1) = (a-b)(a+b-1)$.

If $a,b$ have the same parities then $a+b-1$ is odd and so $2^k \mid a-b$. Since $a,b \in \{0,\dots,2^k-1\}$, in this case $a = b$.

If $a,b$ have different parities then $a-b$ is odd and so $2^k \mid a+b-1$, 
and so $b = 1-a \bmod 2^k$ is uniquely defined, and has a parity different from $a$.
\end{proof}

We can now prove Theorem \ref{thm:main}. First, notice that ${a \choose 2} > a$ for 
$a \geq 4$, and so $a(n) \to \infty$. 
Now suppose that the sequence $\beta$ is ultimately periodic, say with period $\beta(N), \dots, \beta(N+\ell-1)$.
Lemma \ref{lem:bijection} implies that for every $k \ge 1$, the sequence $a(n) \bmod 2^k$ has period 
$a(N) \bmod 2^k, \dots, a(N+\ell-1) \bmod 2^k$, and in particular, 
$a(N) \equiv a(N+\ell) \bmod 2^k$. Choosing $k$ such that $2^k > a(N+\ell)$, we reach a contradiction.

\ifskip\else
\subsection{Normal sequences}
Let $s(n)$ be an integer sequence, and $b \in \N^+$.
The sequence $s^b(n) = s(n) \bmod{b}$ 
is normal, that is, if we chunk it into substrings of length 
$\ell \ge 1$ then each of the $b^\ell$ possible strings of $[b]^\ell$ appear in $s^b(n)$ with equal limiting frequency. 
It is {\em absolutely normal} if it is normal for every $b$.
The sequence $s^b(n) = s(n) \bmod{b}$ can be viewed as a real number $r_b$ written in base $b$.
A classical theorem  from 1922 by E. Borel says that almost all reals are absolutely normal, 
\cite{everest2003recurrence}.
The theorem below shows that MC-finite integer sequences are very rare.
%\cite[xxx]{everest2003recurrence}.

\begin{proposition}
\label{pr:normal}
\begin{enumerate}[(i)]
\item
Almost all reals are absolutely normal. 
\item
If $s^b(n)$ is normal for some $b$, then $s(n)$ is not MC-finite. 
\item
Let $UP_b$ be the set of integer sequences $s^b(n)$ 
with
$s^b(n) = s(n) \bmod{b}$ for some integer sequence $s(n)$
which are ultimately periodic.
Then $UP_b$ has measure $0$.
\end{enumerate}
\end{proposition}
Proving that a specific sequence is normal is usually difficult.

Here is a challenge:
\begin{conjecture}
The binary sequence $\beta(n) = a(n) \bmod 2$ from Theorem \ref{thm:main} is normal with $b=2$. 
\end{conjecture}
%------------------------------------------------------
%\ifskip\else
Here are some possible questions for further research:

\begin{enumerate}[(i)]
\item
Let $p_1 < \cdots < p_\ell$ be prime and $d_1,\dots,d_\ell \ge 0$ be integers. 
Construct a one-dimensional PRS whose reduction modulo $m$ is eventually periodic iff 
$\operatorname{ord}_{p_i}(m) \le d_i$. (Theorem \ref{thm:main} gives an answer for $\ell = 1$, $p_1 = 2$, $d_1 = 0$.)
\item 
\end{enumerate}
\fi %skip

%\bibliographystyle{alpha}
%\bibliography{biblio}
%-------------------------------------------------------------------
%\input{rj-normal}
\subsection{Normal sequences}
Let $s(n)$ be an integer sequence, and $b \in \N^+$.
The sequence $s^b(n) = s(n) \bmod{b}$ 
is normal, that is, if we chunk it into substrings of length 
$\ell \ge 1$ then each of the $b^\ell$ possible strings of $[b]^\ell$ appear in $s^b(n)$ with equal limiting frequency. 
It is {\em absolutely normal} if it is normal for every $b$.
The sequence $s^b(n) = s(n) \bmod{b}$ can be viewed as a real number $r_b$ written in base $b$.
A classical theorem  from 1922 by E. Borel says that almost all reals are absolutely normal, 
\cite{everest2003recurrence}.
The theorem below shows that MC-finite integer sequences are very rare.

Let $PR_b$ be the set of integer sequences $s^b(n)$ 
with $s^b(n) = s(n) \bmod{b}$ for some integer sequence $s(n)$.
$PR_b$ is the projection of all integer sequences to sequences over $\Z_b$.
We think of $PR_b$ as a set of reals with the usual topology and its Lebesgue measure.
Let $UP_b \subseteq PR_b$ be the set of sequences $s^b(n) \in PR_b$ which are ultimately periodic.

\begin{proposition}
\label{pr:normal}
\begin{enumerate}[(i)]
\item
Almost all reals are absolutely normal. 
\item
$s(n)$ is MC-finite iff for every $b \in \N^+$ the sequence $s^b(n)$ is ultimately periodic
\item
If $s^b(n)$ is normal for some $b$, then $s(n)$ is not MC-finite. 
\item
$UP_b \subseteq PR_b$ has measure $0$.
\end{enumerate}
\end{proposition}
Proving that a specific sequence is normal is usually difficult.

Here is a challenge:
\begin{conjecture}
The binary sequence $\beta(n) = a(n) \bmod 2$ from Theorem \ref{thm:main} is normal with $b=2$. 
\end{conjecture}

%-------------------------------------------------------------------
%\input{rj-elim}
%\subsection{Eliminating hard-wired constants}
\section{Eliminating hard-wired constants}
\label{se:constants}

Let $\mathfrak{S}_r(n) =( [r+n], a_1, \ldots , a_r, E)$ be the structures on $[r+n]$ 
where $E$ is an equivalence relation and the $r$ elements $a_1, \ldots , a_r$  are in different equivalence classes.
$S_r(n)$ counts the number of such structures on $[r+n]$.

\ifmargin
\marginpar{s-elim}
\else\fi

Let $\mathfrak{E}_r(n)$ be a structure
on $[n]$ which consists of the following:
\begin{enumerate}[(i)]
\item
$E(x,y)$ is an equivalence relation on $[n]$;
\item
There are $r$ unary relations $U_1, \ldots, U_r$ on  $[n]$;
\item
The sets $U_i(x)$ are disjoint;
\item
Each $U_i(x)$ is either empty or consists of  exactly one equivalence class of $E$;
\end{enumerate}
Let $E_r(n)$ be the number of such structures on $[n]$.

\begin{lemma}
For every $r,n \in \N^+$ there is a bijection $f$ between the structures $\mathfrak{E}_r(n)$ on $[n]$
and the structures $\mathfrak{S}_r(n)$ on $[r+n]$,
hence we have $E_r(n) = S_r(n)$.
\end{lemma}
\begin{proof}
Given a structure $\mathfrak{S}_r(n)$ we define $f(\mathfrak{S}_r(n))$
as follows:
\begin{enumerate}[(i)]
\item
The universe of $f(\mathfrak{S}_r(n))$ is $\{r+1, \ldots, r+n\}$.
\item
If for $i \leq r$ the set $\{i\}$ is a singleton equivalence class, we put $U_i =\emptyset$.
%\item
If there is an equivalence class $E_i$ which strictly contains 
%exactly 
$i$ we put $U_i = E'_i = E_i - \{i\}$.
\item
$E'$ is the equivalence relation induced by $E$ on $\{r+1, \ldots, r+n\}$.
\end{enumerate}
Conversely, given a structure $E_r(n) = ([n], E, U_1, \ldots, U_r)$
we define $g(E_r(n))$ as follows:
\begin{enumerate}[(i)]
\item 
The universe of $g(E_r(n))$ is $[n+r]$ and the equivalence relation $E'$ is defined
by defining its equivalence classes.
\item
If $U_i$ is empty for some $i \geq n+1$ the singleton $\{i\}$  is an equivalence class of $E'$.
%\item
If $U_i$ is not empty, then the equivalence class of $E'$ which contains $i$ is $U_i \cup \{i\}$.
\item
If $C$ is an equivalence class of $E$ such that
$U_i \neq C$ for all $i \geq n+1$,  then $C$ is an equivalence class for $E'$.
\end{enumerate}
It is now easy to check that $f,g$ are bijections and $g$ is the inverse of $f$.
\end{proof}

\begin{remarks}
\begin{enumerate}[(i)]
\item 
Clearly the class of structures $E_r(n)$ as defined here is $\FOL$-definable.
Hence we can apply the Specker-Blatter Theorem and conclude that
$S_r(n)$ is MC-finite.
\item
If $A$ is ultimately periodic
then $S_{A,r}(n)$ is also MC-finite.
To see this we note that
for $S_{A,r}(n)$ all the equivalence classes $C$  satisfy $|C| \in A$.
This means that in a structure $\mathfrak{E}_{A,r}(n)$
the equivalence classes $C$ satisfy $|C| \in A$, if they do not contain a $U_i$,
and $|C| \in A'$ where $A' = \{ a-1: a \in A \}$, otherwise.
If $A$ is ultimately periodic, so is $A'$ and both are definable in $\CMSOL$.
\item
For the Lah numbers $L_{r}(n)$ and $L_{A,r}(n)$ we proceed likewise by replacing the equivalence relation
by a linear quasi-order.
%----------------
For every $i$ we add two  further unary relations and the appropriate conditions 
in order to take care of the ordering of the special elements.
%----------------
\ifskip\else
 you can add for every i two unary relations, one for the x for which a<=a_i and one for the x for which a_i<=x, and add the corresponding restrictions (if I'm not mistaken these are that U_{i,<=}(x) implies U_{i-1,<=(x)}, that  U_{i-1,>=}(x) implies U_{i,<=(x)}, that every x satisfies  U_{i,<=}(x) or U_{i,>=}(x), and that U_{i,>=}(x) and U_{i,<=y}(y) imply the relation between x and y).
\fi %skip
%----------------
Hence both $L_{r}(n)$ and $L_{A,r}(n)$ are MC-finite.
\end{enumerate}
\end{remarks}

%-------------------------------------------------------------------
%\input{rj-c-finite}
\section{Proof of Theorem \ref{th:FM} and its applications}
\label{se:c-finite}
\ifmargin
\marginpar{s-c-finite}
\else\fi

In order to prove Theorem \ref{th:FM} we use Theorem \ref{th:FMapp} below.
For this we have to introduced the definition of $\CMSOL$-definable graph polynomials.

%--------------------------------------------
\ifskip\else
\subsection{Iteratively constructed sequences of graphs}
%\subsection{The Fischer-Makowsky Theorem}
In this subsection we follow \cite{fischer2008linear}.
\begin{definition}
A $k$-colored graph is a graph $G$ together with $k$ sets $V_1,V_2,...,V_k\subseteq V(G)$ such that $V_i \cap V_j =\emptyset $ for $i\neq j$.
A basic operation on $k$-colored graphs is one of the following:
\begin{itemize}
\item $Add_i$: add a new vertex of color $i$ to $G$.
\item $Recolor_{i,j}$: recolor all vertices with color $i$ to color $j$ in $G$.
\item $Uncolor_{i,j}$: remove the color of all vertices with color $i$. Uncolored vertices cannot be recolored again.
\item $AddEdges_{i,j}$: add an edge between every vertex with color $i$ and every vertex with color $j$ in $G$.
\item $DeleteEdges_{i,j}$: delete all edges between vertices with color $i$ and vertices with color $j$ from $G$.
%\item $Duplicate$: Add a disjoint copy of $G$ to $G$.
\end{itemize}
A unary operation $F$ on graphs is $\MSOL$-elementary if $F$ is a finite composition of basic operations 
on $k$ colored graphs (with $k$ fixed). 
We say that a sequence of graphs $\{G_n\}$ is iteratively constructed 
if it can be defined by fixing a graph $G_0$ and defining $G_{n+1}=F(G_n)$ for an $\MSOL$-elementary operation $F$.
\end{definition}

\begin{remark}
In \cite{fischer2008linear} 
there was an additional operation allowed
\begin{itemize}
\item $Duplicate$: Add a disjoint copy of $G$ to $G$,
\end{itemize}
assuming erroneously that $Duplicate$ behaves like a unary operation on graphs.
Although it looks like a unary operation on graphs, the sequence of graphs
$$
G_0 = E_1, G_{n+1} = Duplicate(G_n)
$$
grows too fast and does not fit the framework that the authors have envisaged in \cite{fischer2008linear}.
\end{remark}
\fi %skip
%-------------------------------------------------------------------

\subsection{$\CMSOL$-definable graph polynomials}
\begin{definition}
Let $\Z$ be the ring of integers.  We consider polynomials over 
$\Z[\bar{x}]$. 
For an $\CMSOL$-formula for graphs 
$\phi(\bar{v})$ with $\bar{v}= (v_1, \ldots , v_s)$, define $card_G(\phi)$ to be the cardinality of subsets of $V(G)^s$ 
defined by $\phi$.
The extended $\CMSOL$ graph polynomials are defined recursively. We first
define the {\em extended $\CMSOL$-monomials}.  
Let $\phi(\bar{v}) \in \CMSOL$. 
An extended $\CMSOL$-monomial is a term of one of the following possible forms:
\begin{itemize}
\item $x^{card_G(\phi)}$ where $x$ is one of the variables of $\bar{x}$.
\item $x_{(card_G(\phi))}$ i.e. the falling factorial of $x$.
\item ${x \choose card_G(\phi)}$.
\item $\prod_{\bar{v} \in V(G)^s:\phi(v)}t(\bar{x})$
where $t(\bar{x})$ is a term in $\ZZ[\bar{x}]$.
\end{itemize}
The {\em extended $\CMSOL$ graph polynomials} are obtained from the monomials by
closing  under finite addition and multiplication.
Furthermore they are closed under
summation over subsets of $V(G)$ of the form 
$$\sum_{U:\phi(U)} t$$
where $\phi$ is an $\CMSOL$-formula with free set variables $U$,
and under multiplication over elements of $V(G)^s$ of the form 
$$\prod_{\bar{v} \in V(G)^2:\phi(v)}t(\bar{x})$$.
\end{definition}

\begin{theorem}[Theorem 1 \cite{fischer2008linear}]
\label{th:FMapp}
Let $F$ be an elementary operation on graphs, $\{G_n:n\in \NN\}$ an $F$-iterated sequence of graphs, 
and $P$ an extended $\CMSOL$-definable graph polynomial. 
Then $\{G_n\}$ is C-finite, i.e. there exist polynomials $p_1,p_2,...,p_k\in \ZZ[\bar{x}]$ such that for sufficiently large n,
$$
P(G_{n+k+1})=\sum_{i=1}^k p_iP(G_{n+i})
$$
\end{theorem}

This proves Theorem \ref{th:FM}.

\subsection{Proofs of C-Finiteness}

Now we give the detailed proofs of Theorem \ref{th:c-finite}.
\phantom{fdgf}
%\\
%{\bf Theorem \ref{}} 
\begin{proposition}
\label{R1}
Fix $k_0\in \NN$. Then  $S(n,k_0)$ is a C-finite sequence.
\end{proposition}
\begin{proof}
Let $P$ be the graph property of cliques with at least one vertex, i.e. $P=\{K_n:n\geq 1\}$, and define $G_n=K_n$. 
Note that a $P$-coloring of $G_n$ with $k_0$ colors is a partition of $V(G_n)$ into exactly $k_0$ non empty color classes, 
so $H_P(G_n,k_0)=S(n,k_0)$. We want to apply the Fischer-Makowsky theorem. 
First, note that the sequence $G_n$ is iteratively constructible, 
see Example \ref{exIterativelyConstructed}
or  \cite[Proposition 2]{fischer2008linear}.
Hence $H_P$ is an extended $\CMSOL$ graph polynomial 
and we can use Theorem \ref{th:FMapp}.
\end{proof}

%{\bf Theorem \ref{}} 
\begin{proposition}
\label{R2}
Fix $k_0\in \NN$. Then  $S_r(n,k_0)$ is a C-finite sequence.
\end{proposition}
\begin{proof}
Let $P$ be the graph property of edgeless graphs with at least one vertex, i.e. $P=\{\bar{K_n}:n\geq 1\}$, 
and define $G_n=K_r\cup \bar{K_n}$. Note that a $P$-coloring of $G_n$ with $k_0+r$ colors 
is a partition of $V(G_n)$ into exactly $k_0+r$ non empty color classes, 
such that every vertex in $V(K_r)\subseteq V(G_n)$ is in a different color class, 
so $H_P(G_n,k_0+r)=S_r(n,k_0)$. We want to apply the Fischer-Makowsky theorem. 
First, note that the sequence $G_n$ is iteratively constructible: 
put $G_0=K_r$. Now given $G_n$, we construct $G_{n+1}$ by adding a disjoint vertex.
Hence $H_P$ is again an extended $\CMSOL$ graph polynomial 
and we can use Theorem \ref{th:FMapp}.
\end{proof}

%{\bf Theorem \ref{}} 
\begin{proposition}
\label{R3}
Let $A\subseteq \NN$, and $k_0\in \NN$. Then $S_A(n,k_0)$ is a C-finite sequence if and only if $A$ is ultimately periodic.
\end{proposition}
\begin{proof}
First, note that 
$S_A(n,1)=1$ iff $n\in A$. 
Therefore, if $A$ is not ultimately periodic, 
$S_A(n,1)$ is 
%not ultimately periodic.
not C-finite.
On the other hand, assume $A$ is ultimately periodic. 
Let $P$ be the graph property of cliques with vertex size in $A$, i.e. $P=\{K_n:n\in A\}$, and define $G_n=K_n$. 
Note that a $P$-coloring of $G_n$ with $k_0$ colors is a partition of $V(G_n)$ into exactly $k_0$ non empty 
color classes, with each color class
with size in $A$. 
so $H_P(G_n,k_0)=S_A(n,k_0)$. We want to apply the Fischer-Makowsky theorem. 
As before, the sequence $G_n$ is iteratively constructible.
Hence $H_P$ is again an extended $\CMSOL$ graph polynomial 
and we can use Theorem \ref{th:FMapp}.

\end{proof}

%{\bf Theorem \ref{}} 
\begin{proposition}
\label{R4}
Let $A\subseteq \NN$, and $k_0\in \NN$. 
Then  $S_{A,r}(n,k_0)$ is a C-finite sequence if and only if $A$ is ultimately periodic.
\end{proposition}
\begin{proof}
$S_A(n,1)=1$ iff $n\in A$. 
If $A$ is not ultimately periodic, 
then also $S_A(n,1)$ is 
not C-finite. 
%not ultimately periodic.
Assume $A$ is ultimately periodic. 
Let $P$ be the graph property of edgeless graphs with vertex size in $A$, i.e. $P=\{\bar{K_n}:n\in A\}$, 
and define $G_n=K_r\cup \bar{K_n}$. 
Note that a $P$-coloring of $G_n$ with $k_0+r$ colors is a partition of $V(G_n)$ 
into exactly $k_0+r$ non empty color classes with sizes in $A$, 
such that every vertex in $V(K_r)\subseteq V(G_n)$ is in a different color class, so $H_P(G_n,k_0+r)=S_r(n,k_0)$. 
We want to apply the Fischer-Makowsky theorem. As before, the sequence $G_n$ is iteratively constructible.
Hence $H_P$ is again an extended $\CMSOL$ graph polynomial 
and we can use Theorem \ref{th:FMapp}.
\end{proof}
%-------------------------------------------------------------------
%\input{rj-explicit}
\section{An explicit computation of $S_A(n,k)$}
\label{se:explicit}
\ifmargin
\marginpar{s-explicit}
\else\fi

Let $A\subseteq \NN$. 
$S_A(n,k)$ counts the number of partitions of $[n]$ into $k$ sets with cardinalities in $A$. 
%Problems relating to such partitions (in a more general form) were investigated by 
%Canfield and Wilf in \cite{canfield2012growth}, 
%Ljuji\'c and Nathanson in \cite{ljujic2012partition} and Alon in \cite{alon2013restricted}. 

We shall compute $S_A(n,k)$ explicitly. For $A = \N^+$ this will give also an alternative way of computing
$S(n,k)$, the Stirling numbers of the second kind. 
The method is reminiscent to \cite[Theorem 8.6]{charalambides2018enumerative}
or, in very different notation \cite[Chapter 1, Exercise 45]{bk:Stanley86}.

We introduce some suitable notation.
Let $A\subseteq \NN$. 
$S_A(n,k)$ counts the number of partitions of $[n]$ into $k$ sets with cardinalities in $A$. 
%Let $A\subseteq \NN$ and $a\in \NN$. 
Let $V(A,k)$ be the set of $k$-tuples $(L_1, \ldots, L_k)$
of elements of $A$ ordered in non-decreasing order, 
with $\sum_{i=1}^k =n$,
i.e.
$$
V(A,k)=\{(l_1,l_2,...,l_k)\in A^k:0<l_1\leq l_2\leq...\leq l_k,\sum_{i=1}^k l_i=n\}.
$$
For $(l_1,l_2,...,l_k)\in V(A,k)$ 
define $g(m;l_1,l_2,...,l_k)$ to be the number of times $m$ appears in the $k$-tuple $(l_1,l_2,...,l_k)$, 
and 
$$f(l_1,l_2,...,l_k)=\prod_{m\in(l_1,l_2,...,l_k)}g(m;l_1,l_2,...,l_k)! .$$ 
Next we define inductively:
$c_1 =n$, $c_{i+1} = c_i - l_i$, hence $c_i = n - \sum_1^{i-1} l_i$.

\begin{theorem}
\label{th:explicit}
Let $A\subseteq \NN$.
Then
$$
S_A(n,k)=
\sum_{(l_1,l_2,...,l_k)\in V(A,k)}\frac{1}{f(l_1,l_2,...,l_k)}
\prod_{i=1}^k {c_i \choose l_i}
%{n \choose l_1}{n-l_1 \choose l_2}...{n-l_1-l_2-...-l_{k-1} \choose l_k}
$$
\end{theorem}
\begin{proof}
To partition $[n]$ into $k$ sets with cardinalities in $A$, we proceed as follows: 
First, we select the cardinalities of the $k$ sets. 
This corresponds to picking an element $(l_1,l_2,...,l_k) \in V(A,k)$. 
To construct a partition of $n$, we choose $l_1$ elements from $[n]$, then $l_2$ elements from $[n-l_1]$ etc. 
Finally, we divide by $f(l_1,l_2,...,l_k)$ to account for double counting 
%as the procedure counts ordered partitions.
of tuples with equal entries.
\end{proof}

%-------------------------------------------------------------------
\end{document}